\documentclass{article}

\usepackage[preprint]{neurips_2026}

\usepackage[utf8]{inputenc} 
\usepackage[T1]{fontenc}    
\usepackage{hyperref}       
\usepackage{url}            
\usepackage{booktabs}       
\usepackage{amsfonts}       
\usepackage{nicefrac}       
\usepackage{microtype}      
\usepackage{xcolor}         

\usepackage{multirow}
\usepackage{graphicx}   
\usepackage{caption}
\usepackage{enumitem}

\usepackage{tabularx}

\newcommand{\layoutsNo}{49}
\newcommand{\scenariosNo}{7746}

\newcommand{\set}[1]{\mathcal{#1}}
\newcommand{\norm}[1]{\big\| #1\big\| }

\usepackage[utf8]{inputenc} 
\usepackage[T1]{fontenc}    
\usepackage{hyperref}       
\usepackage{url}            
\usepackage{nicefrac}       
\usepackage{xcolor}         
\usepackage{graphicx}
\usepackage{amsmath,amssymb,amsthm, bm}
\usepackage{dsfont}
\usepackage[english]{babel}
\usepackage{multirow}
\usepackage{longtable}
\usepackage[normalem]{ulem}

\newcommand{\rev}[1]{\textcolor{black}{#1}}

\newcommand{\blue}[1]{\textcolor{black}{#1}}

\usepackage{soul}
\usepackage{appendix}
\usepackage{algorithmic}            
\usepackage{algorithm}
\usepackage{subcaption}

\title{\texttt{WFDroneBench:} A Benchmark for Sensor Placement and Drone Routing for Wildfire Detection}

\author{%
  Romain Puech \\
  Operations Research Center\\
  Massachusetts Institute of Technology\\
  Cambridge, MA, USA \\
  \texttt{puech@mit.edu} \\
  \And
    Joseph Ye \\
  Sloan School of Management \\
  Massachusetts Institute of Technology \\
  Cambridge, MA, USA \\
  \texttt{justdoye@mit.edu} \\
  \And
  Danique de Moor \\
  Sloan School of Management \\
  Massachusetts Institute of Technology \\
  Cambridge, MA, USA \\
  \texttt{demoor@mit.edu} \\
  \And
    Ana Trišović \\
    Computer Science and Artificial Intelligence Lab\\
  Massachusetts Institute of Technology \\
  Cambridge, MA, USA \\
  \texttt{ana\_tris@mit.edu} \\
  \And
  Dimitris Bertsimas \\
  Sloan School of Management \\
  Massachusetts Institute of Technology \\
  Cambridge, MA, USA \\
  \texttt{dbertsim@mit.edu} \\
}

\begin{document}

\maketitle

\begin{abstract}
Increasingly frequent and severe wildfires threaten ecosystems, public health, and infrastructure. Early detection is vital but limited by existing monitoring systems. Drones offer mobile, real-time coverage, but optimizing sensor placement and drone routing in dynamic fire zones remains challenging. To address this, we introduce \texttt{WFDroneBench}, an open-source Python benchmarking library for early wildfire detection that integrates machine-learned risk maps with optimization-based deployment strategies for sensors, charging stations, and drones. It evaluates risk maps, optimization strategies, and monitoring equipment using standardized metrics and realistic wildfire simulations. The framework supports benchmarking across predictive and decision-making components: machine learning researchers can assess risk models \rev{and compare routing strategies}. \texttt{WFDroneBench} includes \scenariosNo{} scenarios across \layoutsNo{} locations, built from historical ignitions, real-world wildfire risk maps, and simulated fire spread, along with two ground detector and \blue{three} drone routing strategies. Our experiments show that \blue{the risk-aware strategy Max-Coverage significantly outperforms other baselines when risk maps
are sufficiently accurate, achieving the fastest detection on the most difficult fires.} We further find that risk-aware static infrastructure helps even under an imperfect
riskmap and drone-based detection outperforms ground sensors. Finally, our results reveal two key open challenges: (i) detecting small fires rapidly and reliably, and (ii) improving risk-map prediction, where the gap between ground-truth ignition patterns and available risk maps highlights a significant opportunity for ML innovation. \blue{We openly release all code, data, and documentation.}
\end{abstract}

\section{Introduction}

Wildfires and the resulting deforestation have been on the rise in recent years~\citep{curtis2018classifying, tyukavina2022global, hu2025coexposure}, causing substantial economic damage and costs exceeding \$70 billion annually~\citep{thomas2017costs}. The 2018 California wildfires alone resulted in \$149 billion in total damages, including \$32 billion in health-related costs from mortality, medical care, and lost productivity~\citep{wang2021economic}. Despite increased wildfire awareness and resource investments, the January 2025 Los Angeles wildfires devastated the wildland–urban interface, causing \$150 billion in damages and evacuating 200,000 residents~\citep{diGiuseppe2025california}. This highlights critical gaps in current detection and response systems. Effective monitoring and early wildfire detection is therefore crucial to ensure fires remain manageable and to minimize damage.

Traditional approaches rely on ground sensors and satellites, but both have notable limitations. Ground sensors provide limited visibility and are costly to build and maintain, while satellites often lack the spatial and temporal resolution needed for timely detection \citep{WildfireDetectionBouguettaya}. Notably, more than half of U.S. counties, particularly in rural regions, lack any air quality monitoring infrastructure \citep{roque2025identifying}. Drones are increasingly recognized as a promising solution for wildfire detection. Their ability to scan large, remote areas with high spatial and temporal resolution makes them well-suited for early fire detection~\citep{WildfireDetectionBouguettaya, WildfireDetectionDLGhali22}. Advances in onboard sensing and deep learning have significantly improved the accuracy of aerial wildfire detection, enabling drones to reliably identify smoke or flames~\cite{WildfireDetectionDLGhali22, WildfireDetectionDLKim, WildfireDetectionDLKumar24, WildfireDetectionDLYandouzi22, WildfireDetectionDLYang}. However, most studies in this domain assume that drones follow pre-defined flight paths. The critical question of how to actively plan drone routes to maximize early discovery of fires remains largely unaddressed. 

Most existing work on drone routing focuses on post-detection monitoring and response 
(see e.g., \cite{DroneroutingPostGhamry17, WildfireMonitoringBailonruiz22, DroneroutingPostGhamry17}). 
In contrast, proactive wildfire detection, where drones search for unknown fires across large landscapes, is far less explored. Existing detection efforts largely center on the placement of ground sensor networks \cite{BAO2015100, f12040453, 10756651, 9068226} or on ML-driven image analysis to enhance onboard fire recognition by drones (see e.g., \cite{RAMADAN2024101248, WildfireDetectionBouguettaya}), neglecting the integration of active routing strategies. Yet, successful early detection demands tight coupling between intelligent routing and accurate sensing. This includes integrating mobile drones and stationary sensors, dynamically adapting to changing risk conditions, and operating under constraints like battery life, transmission range, and real-time data availability. Key challenges include: how to route drones to maximize early detection probability, how to adapt to evolving risk, and how to coordinate drone and sensor deployments under limited resources.

In this paper, we address the underexplored challenge of proactive wildfire detection. We introduce an integrated framework for optimizing and evaluating risk-aware sensor placement and drone routing, supported by a modular open-source Python library with standardized evaluation metrics and semi-synthetic benchmarking datasets. Our main contributions are:



\begin{enumerate}
    \item \textbf{Introducing a problem formalism for wildfire detection with ground sensors and aerial drones.} We develop a mathematically-grounded framework for proactive wildfire surveillance that integrates evolving risk information, and real-world constraints such as drone battery life, recharging, transmission range, and multi-drone coordination. This enables designing adaptive, risk-aware routing strategies designed to improve the timeliness and effectiveness of wildfire detection efforts. 
    \item \textbf{\texttt{WFDroneBench}, a benchmarking framework for wildfire detection strategies.} \texttt{WFDroneBench} is a flexible, open-source framework that integrates wildfire risk data, drone profiles, and operational constraints. It supports the development and evaluation of detection strategies, risk maps, and drone specifications using diverse metrics, enabling controlled ablations and cross-method comparisons. The framework facilitates experiment reproducibility, visual diagnostics, and community-driven development, and supports benchmarking of both predictive and decision-making components across varied use cases.

    \item \textbf{A realistically diverse, semi-synthetic wildfire benchmark dataset.} We curate a semi-synthetic dataset combining historical ignition records with physics-based wildfire spread simulations to represent fires spreading without human intervention. Scenarios span diverse geographic and ecological conditions, capturing variations in climate and weather conditions, terrain, ignition, wildfire, size and rate of propagation, all captured in rich metadata enabling detailed subgroup analysis for robust benchmarking.  
    \item \textbf{Developing, implementing, and benchmarking diverse wildfire detection strategies in a uniform code base.} We adapt, implement and evaluate \blue{six} wildfire monitoring strategies: two for ground sensor placement and \blue{three} for drone routing. These include both heuristic and optimization-based methods, all implemented within a consistent Julia codebase. Our experiments reveal that high-accuracy risk maps are essential for enabling effective risk-aware strategies, risk-aware static infrastructure helps even under an imperfect riskmap, and drone-based detection outperforms ground sensors.
\end{enumerate}




\section{Related Work}

Previous work on optimal drone routing for wildfire detection spans a range of modeling approaches. One line of work proposes mixed-integer linear programs (MILPs) to jointly optimize the placement of watchtowers or balloons alongside battery-constrained, single-route drone missions~\citep{DroneroutingDeLaFuente}. While these methods note that drones can be reused, recharging is not explicitly modeled. The probabilistic path planning (PPP) strategy proposed in \cite{fire7070254} uses logistic regression to estimate fire probabilities and subsequently employs a dynamic programming (DP) algorithm to generate an optimal routing for unmanned aerial vehicles (UAV). However, it does not explicitly consider operational constraints such as battery life, charging requirements, or return-to-base logistics. Other studies minimize the makespan, defined as the longest mission duration across all drones, motivated by constraints such as connectivity synchronization, coordinated data retrieval, and bandwidth limits; in some cases, this is coupled with the placement of recharging stations \citep{DroneroutingSantin, DroneroutingJemmali}. These approaches are well-suited for structured surveillance or exhaustive coverage, where all regions are considered equally important. An alternative strategy, the risk-specific UAV patrol path (RSUPP), uses Gaussian mixture clustering to partition high-risk regions into subareas, assigning each UAV to the shortest route through one such region based on static wildfire risk~\citep{DroneroutingXu}.
Other modeling approaches include nature-inspired metaheuristics such as Particle Swarm Optimization (PSO) \cite{DroneroutingYan}. There are many approaches for generic UAV routing (see e.g. \citep{drones9040263, s23031463, rs16214019}).

We implement \blue{three} proactive drone routing approaches that leverage wildfire risk data, accounting for spatial heterogeneity and enabling adaptive replanning as conditions evolve. Unlike models assuming uniform importance or known fire locations, the incorporated methods adapt to evolving wildfire risk. 

\section{Problem Formulation}
\label{sec:formulation}

Our objective is to effectively deploy a limited number of ground sensors, charging stations, and drones across a real-world geographic region for early wildfire detection. To formalize the problem, we consider a spatiotemporal environment represented as a two-dimensional grid $\mathcal{I}$ of size $N \times M$ over a discrete time horizon $\set T = \{1, \dots, T\}$. At each time step $t \in \set T$, a fire may ignite and propagate to adjacent grid points. If a fire ignites at grid point $i$ and time $t_0$, it is considered detected once it enters the sensing radius of any deployed device at time $t_\text{det}$, resulting in a detection delay $\Delta t = t_\text{det} - t_0$.

To minimise the expected detection delay $\mathbb{E}[\Delta t]$ over the ignition distribution, we define a \textbf{detection strategy} that optimizes the placement of ground sensors and charging stations, denoted by $\{x_i^g, x_i^c\}$, and the routing of and charging of drones over time, denoted by $\{a_{its}, c_{its}\}$. This strategy is guided by a spatiotemporal wildfire risk map, which informs both infrastructure deployment and drone movement decisions. The fire risk is encoded as a burn risk field $r_{it} \in [0,1]$, representing the probability that grid point $i$ is burning at time $t$, incorporating both ignition and spread probability. This risk map can be obtained from physical simulation models or predicted using machine learning methods that capture the underlying dynamics of wildfire behavior.

Devices may only be deployed at feasible subsets: ground sensors on $\set I_g \subseteq \set I$ and charging stations on $\set I_c \subseteq \set I$. We define binary placement variables $x_i^g \in \{0,1\}$ and $x_i^c \in \{0,1\}$ to indicate whether a sensor or station is placed at location $i$. A drone fleet $\set S$ is routed through the domain according to binary variables $a_{its}$ and $c_{its}$, where $a_{its} = 1$ if drone $s$ flies at location $i$ at time $t$ and $c_{its} = 1$ if drone $s$ is charging at grid point $i$ at time $t$ (see notation overview in Table~\ref{tab:sets}.) Device allocation is constrained by resource limits. Spatial exclusion constraints further prevent deployments from being placed too close to one another. Drones follow a single-hop kinematic rule—after spatial and temporal rescaling using the drone cruising speed, they may move only to neighboring cells at each time step—and must remain within a communication radius $R$ of a charging station. Each drone’s endurance is limited to $D_{\max}$ time steps and is fully replenished upon landing at a charging station.

\begin{table}[!htb]
\centering
\resizebox{\textwidth}{!}{%
\begin{tabular}{ll|ll} \toprule
\textbf{Sets}                  & \textbf{Description}             & \textbf{Variables} & \textbf{Description}          \\ \midrule
$\set I$ & Set of all possible grid points on the map & $x_i^g \in \{0,1\}$ & 1 if a ground sensor is placed at grid point $i$ \\
$\set I_g$ & Set of all feasible grid points for the ground sensors & $x_i^c \in \{0,1\}$ & 1 if a charging station is placed at grid point $i$ \\
$\set I_c$ & Set of all feasible grid points for the charging stations & $a_{its} \in \{0,1\}$ & 1 if drone $s$ flies at grid point $i$ at time $t$  \\
$\set I_d$ & Set of all feasible grid points for the drones & $c_{its} \in \{0,1\}$ & 1 if drone $s$ is charging at grid point $i$ at time $t$ \\
$\set T$ & Total time horizon with cardinality $T$ & $\theta_{it} \in \{0,1\}$ & 1 if grid point $i$ is covered by a drone at time $t$ \\
$\set H$ & Optimization time horizon with cardinality $H$ & $b_{ts} \in \mathbb{Z}$ & Battery of drone $s$ at time $t$ (unit time) \\
$\set S$ & Set of drones &   \\ \cmidrule{1-2}
\textbf{Parameters}                  & \textbf{Description}             &           \\ \cmidrule{1-2}
$r_{it}^d$ & Dynamic wildfire risk in grid point $i$ at time $t$ &    \\
$r_i^s$ & Static wildfire risk in grid point $i$ \\
$T_{\text{reev}}$ & Reevaluation step \\
\bottomrule \addlinespace
\end{tabular}}
\caption{Notation used in the problem formalism and the formulation in Section~\ref{sec:strategies}.}
\label{tab:sets}
\end{table}
\section{\texttt{WFDroneBench}}

\begin{figure}[!t]
     \centering              \includegraphics[width=.85\textwidth]{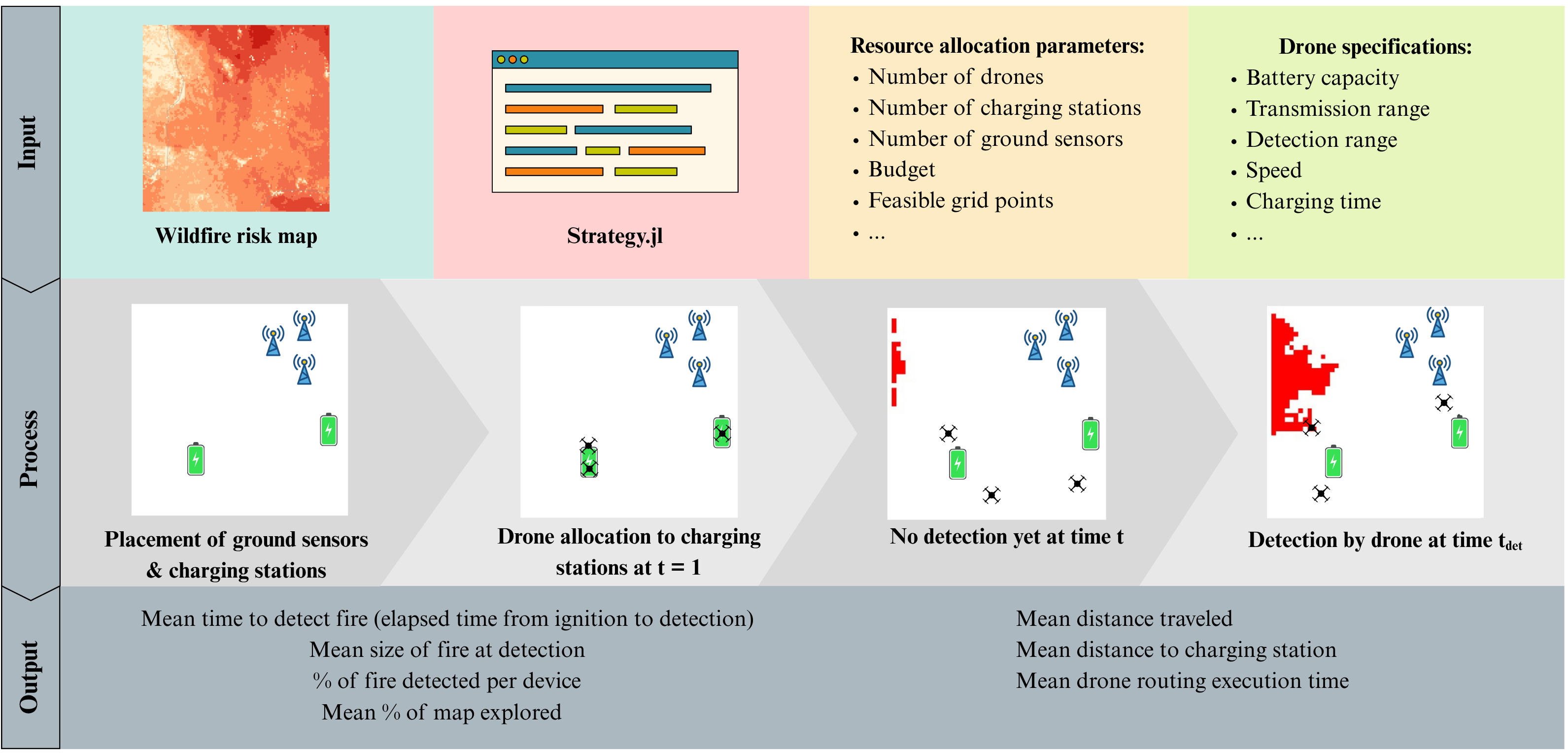}
         \caption{\blue{Workflow of the integrated wildfire monitoring optimization framework.
Given a wildfire risk map, a strategy, and system parameters (number of drones, charging stations, sensors, and drone specifications), the model optimizes (i) infrastructure placement and (ii) drone allocation and routing. The framework simulates detection over time and evaluates performance via detection time, coverage, travel distance, and computational efficiency metrics.}}

         \label{fig:optimizationframework}
\end{figure}

\paragraph{Design goal.} To bridge data-driven risk modeling with operational decision-making, we introduce \texttt{WFDroneBench}, an open-source benchmark for early wildfire detection. The framework supports the integration of machine learning–generated risk maps with optimization-based deployment strategies for ground sensors, charging stations, and drones (see Figure~\ref{fig:optimizationframework}). It includes tools for loading wildfire risk maps, simulating dynamic wildfire scenarios, implementing deployment strategies, and evaluating early detection performance through standardized metrics. \texttt{WFDroneBench} is modular and extensible. Users can implement custom deployment policies or integrate their own risk maps and fire scenarios to evaluate end-to-end performance. The framework supports benchmarking across both predictive and decision-making components across diverse use cases: ML researchers can benchmark risk models via their impact on detection with fixed routing and OR experts can compare routing strategies on a shared risk map.

\paragraph{Framework overview.}\label{terminology} 
The benchmarking dataset is organized into \textbf{layouts}, each representing a distinct geographic map. Each layout is associated with multiple wildfire \textbf{scenarios}, defined as spatiotemporal grids $\set I \times \set T$ that capture wildfire spread over time along with environmental conditions such as weather and vegetation. Scenarios provide the spatial and temporal context in which detection strategies—covering both device placement and drone routing—are evaluated.

To ensure consistency between simulation dynamics and decision-making, we align each scenario’s temporal structure with drones’ operational capabilities by distinguishing between two levels of resolution: \textbf{data} and \textbf{operational}. The \emph{data temporal resolution} defines the duration of each simulation step (e.g., one hour), while the \emph{data spatial resolution} specifies the physical size of each grid point in $\set I$. In contrast, the \textbf{operational resolution} defines the granularity of decisions. The \textbf{operational spatial resolution} corresponds to the distance unit for drone routing, typically a multiple of the drone’s coverage radius. The \textbf{operational time resolution} indicates how many operational grid cells a drone can traverse per simulation step.

A \textbf{strategy} refers to a decision-making policy for wildfire detection, comprising sensor placement and drone routing. The sensor placement component determines where ground sensors and charging stations are deployed, while the drone routing component governs the movement of drones across the grid. Each strategy operates on a \textbf{wildfire risk map} where $r_{it} \in [0,1]$ represents the burn probability of grid point $i$ at time $t$. Our library distinguishes between two types of risk maps: \emph{static} (fixed over time) and \emph{dynamic} (continuously updated in response to changing weather or environmental conditions). During runtime, the risk map is continuously adjusted, lowering ignition and burn probabilities in regions currently monitored by drones to model the mitigating impact of real-time aerial surveillance.

Drones operate under physical and operational constraints. These include finite battery life ($D_{\max}$), limited transmission range ($R$), and bounded mobility. The \textbf{maximum effective speed} of a drone does not reflect raw flight velocity but instead the rate at which it can effectively scan and update grid points within its coverage area. These constraints, together with the scenario resolution, define the feasible range of strategic decisions.

\subsection{Benchmark Datasets}\label{subsec:benchmarkdatasets} 


\paragraph{Data collection.} To model wildfire spread, we use the Sim2Real-Fire dataset~\citep{NEURIPS2024_02e978a2}, which contains over one million high-fidelity simulated fire trajectories and 1,000 real-world wildfire cases, each annotated with five aligned modalities: topography, vegetation, fuel, weather, and satellite imagery. While rich in fire dynamics, its ignition points are synthetically and uniformly distributed, limiting realism. To address this, we incorporate historical ignition data from the Fire Program Analysis Fire-Occurrence Database (FPA FOD)~\citep{short2021spatial}, which includes 2.17 million geo-referenced ignition records across the U.S. This integration enables more realistic study of ignition and spread behavior while avoiding confounding effects from suppression, using counterfactual simulations from physics-based and empirical models. We constrain matches to ignition points within 150 meters of historical locations and allow any calendar day, further limiting selection to large-scale regions at least 15 kilometers wide for operational relevance. The final curated dataset includes \scenariosNo{} fire-spread scenarios across \layoutsNo{} rectangular layouts, with final dimensions ranging from 30 km to 60 km.



\paragraph{Benchmarking dataset categorization.}  
To enable structured analysis of drone deployment strategies under diverse wildfire conditions, we categorize each fire-spread scenario along five key axes: (i) fire size (small vs. large, with large defined as a final burn radius exceeding 20 km), (ii) spread rate (slow vs. fast, where fast denotes rapid expansion within the first 10 hours), (iii) historical match (alignment in both spatial location and calendar day with recorded wildfire events), (iv) seasonal match (correct region and seasonality), and (v) layout inclusion (all fires occurring within a given deployment layout). This categorization facilitates the evaluation of risk maps and response strategies across a range of fire behaviors and temporal conditions, and is especially important for assessing robustness to both localized and large-scale events (e.g., 20–50 km fires)\citep{WildfireriskXu2024}. To avoid temporal bias from fixed ignition times, we introduce a random delay of one to twelve hours in each scenario. This ensures that drone systems begin in a stable operational state and prevents routing algorithms from exploiting prior knowledge of fire onset time.




\subsection{Evaluation Metrics}\label{subsec:evaluationmetrics}



To evaluate monitoring strategies, we measure detection time, coverage, and computational cost. The main detection metric is the ~\textbf{mean detection time} or the expected detection delay $\mathbb{E}[\Delta t] = \frac{1}{|\mathcal{W}|} \sum_{(i, t_0) \in \mathcal{W}} (t_d(i) - t_0(i))$, averaged over wildfire scenarios $\mathcal{W}$ per each layout. To focus on early detection, we consider a fire undetected if it is not spotted within its first twelve hours after ignition. We report the percentage of  detected fires, the mean fire size at detection, and the proportion of detections attributed to each device type (ground sensors, charging stations, and drones). Spatial efficiency and resource metrics include the mean percentage of the map explored, mean distance traveled, and mean distance to the nearest charging station. To account for algorithmic efficiency, we report the mean drone routing execution time, measuring computational overhead. These metrics provide a comprehensive basis for evaluating trade-offs between responsiveness, resource utilization, and computational cost (see Appendix \ref{appendix:evaluation metrics} for full definitions of evaluation metrics).

\begin{figure}
    \centering
\includegraphics[width=\textwidth]{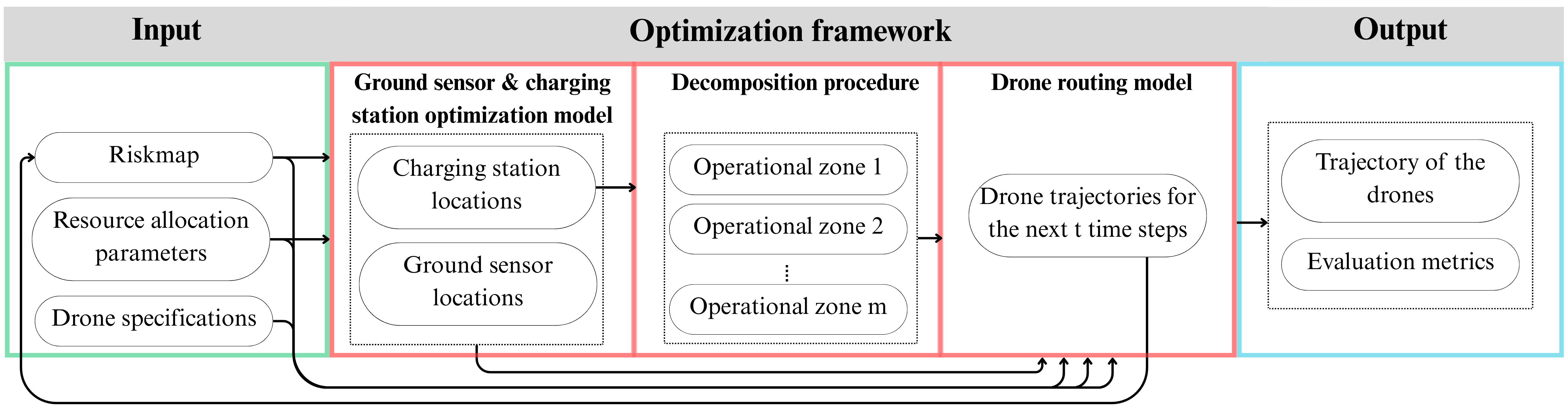}
    \caption{Architectural diagram showing the benchmark components and how they relate to one another.}
    \label{fig:diagram}
\end{figure}

\section{Experiments} ~\label{sec:experiments}

We conducted a comprehensive set of experiments to evaluate wildfire detection strategies under realistic, operational constraints. Our aim was to assess the performance of various sensor placement and drone routing methods across diverse terrains and ignition patterns. This section outlines the detection strategies and baselines, wildfire risk map selection, drone specs, and experimental setup, including hardware configurations.

\subsection{Wildfire Detection Strategies and  Baselines}\label{sec:strategies}
We implement several wildfire detection strategies, consisting of the placement of ground sensors and charging stations, and the routing of drones, and systematically evaluate their performance. To ensure scalability, charging station locations serve as anchors for a decomposition algorithm (see Appendix \ref{appendix:decomposition alg}) that partitions the surveillance area into operational zones, each assigned to a subset of drones. An overview of the framework is shown in Figure~\ref{fig:diagram}.

\paragraph{Ground Sensor and Charging Station Placement Baselines}

\begin{enumerate}[leftmargin=1.1em, itemsep=0.3em, topsep=0.5em]\label{gaussian-coverage}
    \item \textbf{Gaussian-weighted Max-Coverage} (GaussianCov): We model the placement of ground sensors and drone charging stations as a maximum risk coverage location problem, first introduced by~\cite{MCLPChurchRevelle}. We extend the objective to maximize the total monitored wildfire risk, i.e., $$\max_{(\bm{x^g}, \bm{x^c}) \in \mathcal{F}_1} \sum_{t \leq H}\sum_{i \in \mathcal{I}} r^d_{it} \dot \min(1, K_i(x_i^c), x^g_i)$$ with $\mathcal{F}_1 = \{ (\bm{x^g}, \bm{x^c}) \in \{0,1\}^{\mathcal{I}} \times \{0,1\}^{\mathcal{I}} \ \vert \ x^g_i = 0 \ \forall i \in \mathcal{I} \setminus \mathcal{I}_g \text{ and }  x^c_i = 0 \ \forall i \in \mathcal{I} \setminus \mathcal{I}_c    \}$ the set of all valid decisions, $r^d_{it}$ the wildfire risk at cell $i$ and time $t$, and $H$ the optimization horizon. Our approach extends traditional max coverage by introducing a kernel function $K_i$ that captures both direct and indirect coverage effects. For ground sensors we consider standard binary coverage of their immediate location. For charging stations, however, we account for the indirect impact on risk coverage through the drones that will utilize them. Beyond covering their immediate cell, charging stations enable drone operations that can monitor surrounding areas. To model this accessibility, we compute $K_i(x_i^c)$ as the probability that a drone starting from charging station locations will visit cell $i$, based on a 2D Brownian motion model. This yields a Gaussian distribution centered at each station, with variance determined by drone battery life. Overlapping coverage is handled by summing individual probabilities, capped at 1. We estimate $K_i$ via iterative $3 \times 3$ convolutions with a uniform kernel to a grid initialized with zeros and ones to simulate diffusion (see Appendix~\ref{appendix:max coverage ground} for details).

    \item \textbf{Uniform} (Random): We draw uniformly at random the location of each ground sensor and charging station across the grid. This approach ignores wildfire risk, spatial structure, and feasibility constraints. While impractical for real-world use, it provides a useful baseline to assess the benefits of informed, risk-aware placement strategies.
\end{enumerate}

\paragraph{Drone Routing Baselines}
\begin{enumerate}[leftmargin=1.1em, itemsep=0.3em, topsep=0.5em]
    \item \textbf{Max-Coverage} (MaxCov): 
    We formulate drone routing as a rolling-horizon, risk-weighted maximum unique coverage optimization problem. The objective is to maximize cumulative risk observed across all drones and time steps within a planning horizon: $\max_{(\bm a, \bm b, \bm c, \bm \theta) \in \set F_2} \sum_{i \in \set I_d} \left\{r_{i1}^d \theta_{i1} + \sum_{t \in \set H \setminus \{1\}} r_{it}^d (\theta_{it} - \theta_{i,t-1}) \right\}, $ where $\set F_2$ represents the set of feasible decisions subject to drone dynamics and operational constraints (formally defined in Appendix \ref{appendix:max coverage drones}). Wildfire risk scores $r_{it}^d$ are updated every $T_{\text{reev}}$ steps based on drone observations and external data sources. This strategy plans over a finite optimization horizon $\set H$ to identify routing decisions that maximize first-time surveillance of high-risk areas, while accounting for long-term feasibility, including battery depletion, and station returns. 


    \item \textbf{Uniform Coverage} (UniCov):
    As an "uninformed" stochastic baseline, drones are directed to locations sampled uniformly at random within their assigned operational zone. We consider a drone routing strategy for maximum coverage under a uniform risk map, where every grid cell is equally important. This approach, which ignores wildfire risk and prior coverage, serves as a baseline to quantify the added value of risk maps in drone routing.
    \item \textbf{Brownian} (Brownian): At each step, the next waypoint is sampled with probability proportional to the current risk level. Drones are initially placed uniformly across charging stations and move randomly to neighboring grid points at each time step. All movements are constrained to ensure feasibility with respect to battery levels and charging station accessibility. This strategy ignores wildfire risk and provides a simple, constraint-satisfying baseline.
\end{enumerate}

\subsection{Wildfire Risk Map Specification}~\label{sec:risk_maps}

In our experimental setup, we incorporate a real-world, static wildfire Burn Probability (BP) risk map developed by the USDA Forest Service for all U.S. lands~\citep{dillon2023whp}. BP represents the annual probability of wildfire occurrence at a specific location. The dataset provides high spatial resolution at 30 meters, enabling fine-grained, risk-aware planning for both ground sensor placement and drone routing. It is widely adopted in federal and community-level wildfire risk assessments and is publicly available alongside comprehensive methodological documentation.

Our library also supports dynamic wildfire risk maps that evolve over time. To demonstrate this capability and establish a benchmark for evaluation, we compute a mock “ground-truth” burn map using a frequentist approach. For each layout, we average the simulated fire propagation scenarios over time to estimate the empirical burn probability of each cell at each time step, corresponding to $r^d_{i,t}$ in our framework. This dynamic map enables evaluation of detection strategies under ideal conditions and allows comparison with the static BP map to assess sensitivity to risk map accuracy.

\subsection{Drone Specification}

We model wildfire monitoring using real-world drone specifications to ensure operational feasibility. Our experiments use the Mugin EV350, a long-range, fully electric VTOL drone designed for autonomous environmental monitoring. It offers a 130 km flight range, 120-minute flight time, and a cruising speed of 24 m/s, making it suitable for rapid response over large areas. It supports a 2 kg payload, has a transmission range of up to 80 km, and withstands wind conditions up to Grade 6~\cite{mugin_ev350}.

\subsection{Experimental Setup}

To evaluate the effectiveness of our sensor placement and drone routing strategies under realistic wildfire conditions, we design a suite of experiments grounded in operational constraints and historical fire data. Each scenario deploys \textbf{eight} ground stations, \textbf{two} drones, and \textbf{two} charging stations to simulate a realistic and resource-constrained wildfire monitoring operation. To ensure consistency and comparability, all experiments use standardized specifications for the drone platform, charging infrastructure, and ground sensors (see Appendix~\ref{appendix:sets} for full list of variables and Appendix~\ref{appendix:par values} for their specific values). For testing, we use a subset of our dataset (introduced in Section~\ref{subsec:benchmarkdatasets}), which has 80\% or more historical ignition and simulation (Sim2Real) matches, resulting in 12 layouts with 474 wildfire spread scenarios in total. We adopt a high match threshold of 80\% to improve the alignment with real-world risk patterns and to maintain statistical consistency across scenarios.

Our strategies (introduced in Section~\ref{sec:strategies}) comply with operational constraints of planning drone routing in one-hour increments. The runtime remains efficient across methods. For \blue{all strategies}, each hour of operation requires less than 30 seconds of computation. 
Ground station placement strategies require less than two seconds.

\paragraph{Hardware specification.} All our experiments are conducted in parallel on a cluster of Intel Xeon Platinum 8260. Each burn map and sensor placement strategy combination is run as a single job on 32 CPUs with 125 GB of RAM. Within each job, the drone strategies are run sequentially. For each drone strategy, every layout is run in parallel on a separate thread. The routing and placement strategies are implemented in Julia version 1.11.3 using the \texttt{JuMP.jl} version 1.26.0 optimization modeling package and solved with Gurobi version 12.0.3 \citep{Gurobi}.

\section{Results and Discussion} \label{sec:results}

Our experimental results using a real-world wildfire BP map are shown in Table~\ref{tab:bp_results}. To illustrate the outcomes under an ideal risk map, assess sensitivity to risk map quality, and demonstrate the use of the dynamic risk map available in our library, we present results using our “ground truth” risk map (described in Section~\ref{sec:risk_maps}) in Table~\ref{tab:ground_truth}. This map serves as a high-fidelity benchmark with hourly resolution, enabling us to evaluate the upper bound of performance for risk-informed strategies. These results reveal several key insights.

Most notably, the best strategy depends on risk-map fidelity -- \blue{Max-Coverage} excel when risk maps are accurate, whereas spread-out coverage performs better under a less informative risk map.
From Table~\ref{tab:bp_results}, we observe that Uniform Coverage frequently yields the highest detection rates -- and even achieves perfect detection (100$\%$) in all fast-fire scenarios -- while maintaining reasonable detection times. By assuming all areas have the same risk, Uniform Coverage avoids the risk of missing fires due to inaccurate risk assessments. In contrast, Max-Coverage can detect a subset of fires very quickly (e.g., 2.00 hours detection in Fast Big) when the risk map is accurate, but misses many fires when the BP risk prediction provides unreliable guidance (e.g., 73$\%$ versus $100\%$ detection in Fast Big). 
Importantly, Table~\ref{tab:ground_truth} shows the opposite behavior under the ground-truth map: when risk information is perfect, Max-Coverage outperforms Uniform Coverage on detection speed while performing similarly on overall detection rate. This demonstrates that risk-aware routing can be substantially more efficient when risk maps are sufficiently accurate, as Uniform coverage cannot exploit the efficiency gains that come from focusing on high-risk regions. The dependence on risk-map fidelity highlights a central challenge in wildfire monitoring: optimizing drone routing requires matching the routing strategy to the quality of risk information available.



Second,
the comparison between Max-Coverage and Uniform Coverage further highlights the importance of the risk map. Since Uniform Coverage is an ablation of Max-Coverage that follows the same routing logic but assumes a uniform risk distribution, differences in their performance directly reflect how effectively Max-Coverage leverages the input risk information. When using the dynamic ground-truth map, Max-Coverage consistently outperforms Uniform Coverage in terms of detection time. However, this advantage largely disappears when the static BP risk map is used. This trend demonstrates the critical importance of accurate, temporally dynamic risk maps in supporting the effective operation of risk-informed strategies.

Third, risk-aware static infrastructure helps even under an imperfect risk map, and becomes even more beneficial when the map is accurate. This can be observed from the consistent improvement of GaussianCov over Random sensor placement in both Table~\ref{tab:bp_results} and \ref{tab:ground_truth} as in nearly all fire types, GaussianCov yields higher detection rates and often faster detection times because the ground sensors and charging stations are concentrated around regions identified as high risk.

Fourth, drone-based detection consistently outperforms ground sensors, with drones identifying most fires despite only two drones versus eight sensors and two charging stations (Appendix~\ref{appendix:evaluation metrics results}). This suggests drones may offer better returns on investment due to their greater coverage. 

Finally, unsurprisingly, the Brownian strategy is the worst-performing baseline. \blue{Additional data plots are provided in Appendix \ref{appendix:dataset_plots}.}

\renewcommand{\arraystretch}{0.9}

\begin{table}[htbp]
\centering
\tiny

\begin{tabular}{@{}llp{1.5cm}p{1.5cm}p{1.5cm}p{2cm}p{2cm}@{}}
\toprule
 & & Overall & Slow Small & Slow Big & Fast Small & Fast Big \\
\midrule
\textbf{Sensor Strategy} & \textbf{Drone Strategy} & & & & & \\
\midrule

\multirow{4}{*}{Random} 
  & MaxCov   & 91\%; 3.19 (4.42) & 78\%; 3.75 (3.40) & 93\%; 3.54 (4.40) & \textbf{100\%}; 2.29 (3.40) & 55\%; 9.67 (9.63) \\
  & UniCov   & \textbf{93\%}; 3.08 (4.23) & \textbf{88\%}; 3.97 (4.84) & \textbf{94\%}; \textbf{3.21} (4.12) & \textbf{100\%}; 3.00 (2.08) & \textbf{100\%}; \textbf{7.45} (5.35) \\
  & Brownian & 74\%; 5.22 (6.71) & 46\%; 5.37 (6.99) & 77\%; 7.66 (6.52) & 57\%; 5.25 (4.72) & 82\%; 8.78 (6.28) \\
\midrule

\multirow{4}{*}{GaussianCov} 
  & MaxCov   & 96\%; \textbf{3.16} (4.24) & 90\%; 5.11 (5.65) & \textbf{97\%}; 3.23 (4.08) & \textbf{100\%}; 3.71 (6.90) & 73\%; \textbf{2.00} (2.14) \\
  & UniCov   & \textbf{97\%}; 3.17 (4.17) & \textbf{93\%}; 5.05 (5.45) & \textbf{97\%}; \textbf{3.18} (3.97) & \textbf{100\%}; \textbf{1.00} (1.00) & \textbf{100\%}; 5.73 (6.17) \\
  & Brownian & 81\%; 4.99 (6.33) & 59\%; 6.96 (7.42) & 83\%; 6.34 (6.31) & 86\%; 8.00 (4.69) & 91\%; 6.10 (3.51) \\

\bottomrule
\end{tabular}
\caption{Percentage of fires detected in under 12 hours and mean detection time (±SD in hours) using the BP risk map, across sensor/drone strategies and fire scenarios (by size and speed). }\label{tab:bp_results}
\end{table}

\begin{table}[htbp]
\centering
\tiny
\begin{tabular}{@{}llp{1.5cm}p{1.5cm}p{1.5cm}p{2cm}p{2cm}@{}}
\toprule
 & & Overall & Slow Small & Slow Big & Fast Small & Fast Big \\
\midrule
\textbf{Sensor Strategy} & \textbf{Drone Strategy} & & & & & \\
\midrule

\multirow{4}{*}{Random} 
  & MaxCov   & 89\%; \textbf{2.49} (3.74) & 78\%; 3.72 (4.39) & 89\%; \textbf{2.88} (3.66) & \textbf{100\%}; 2.14 (3.67) & \textbf{100\%}; 4.55 (5.24) \\
  & UniCov   & \textbf{92\%}; 3.74 (5.11) & \textbf{93\%}; 4.71 (5.04) & 92\%; 4.19 (5.22) & \textbf{100\%}; \textbf{1.86} (4.06) & 91\%; 4.20 (4.05) \\
  & Brownian & 66\%; 3.77 (5.55) & 29\%; 7.83 (4.73) & 68\%; 6.16 (5.38) & 43\%; 7.67 (1.15) & 91\%; 5.40 (5.93) \\
\midrule

\multirow{4}{*}{GaussianCov} 
  & MaxCov   & 95\%; 3.20 (4.21) & 78\%; 4.44 (6.16) & 96\%; \textbf{3.32} (4.01) & \textbf{100\%}; 3.57 (3.21) & \textbf{100\%}; \textbf{5.73} (5.50) \\
  & UniCov   & \textbf{96\%}; 3.41 (4.39) & \textbf{83\%}; 3.35 (3.60) & \textbf{97\%}; 3.58 (4.40) & \textbf{100\%}; 2.14 (2.48) & 82\%; 8.22 (5.93) \\
  & Brownian & 73\%; 5.11 (6.53) & 63\%; 7.92 (5.93) & 75\%; 7.23 (6.34) & 86\%; 9.00 (7.13) & 36\%; 8.75 (5.32) \\

\bottomrule
\end{tabular}
\caption{Percentage of fires detected in under 12 hours and mean detection time (±SD in hours) using the ``ground truth'' dynamic map, across sensor/drone strategies and fire scenarios (by size and speed). }\label{tab:ground_truth}
\end{table}

\section{Conclusion}\label{sec:conclusion}

\blue{\texttt{WFDroneBench} is a flexible, extensible, open-source benchmarking framework for early wildfire detection that jointly optimizes ground sensor placement and drone routing under realistic operational constraints. It enables reproducible evaluation across multiple dimensions, including wildfire risk modeling (e.g., dynamic, learning-based risk maps), deployment strategies, and operational decision-making. The framework also incorporates drone-specific characteristics such as battery capacity, transmission range, and detection coverage, bridging algorithmic development with real-world deployment considerations. It fosters interdisciplinary collaboration across machine learning, operations research, and environmental science, and lays the foundation for engagement with public agencies and field responders working toward scalable, data-driven wildfire resilience.}

\blue{Our results show that risk-aware routing (Max-Coverage) significantly outperforms baselines when accurate dynamic risk maps are available, and that drone-based detection is generally more effective than ground sensors. We identify that the bottleneck for drone-based wildfire detection systems lies in developing accurate, temporally dynamic burn probability maps rather than routing optimization alone, calling for more research in Machine Learning devoted to the creation of high-resolution risk maps.}


\paragraph{Limitations.} Despite its strengths, \texttt{WFDroneBench} has several limitations. The simulated environment abstracts away operational complexities such as communication failures, terrain constraints, and execution uncertainties. Additionally, the current implementation assumes full observability and reliable risk updates, which may not hold in real deployments. 
Finally, environmental factors (e.g., wind, temperature) that may directly impact drone performance are not considered. Since such factors typically vary little within an hour, they can be incorporated into the rolling-horizon framework without jeopardizing the drones’ ability to return to a charging station, as long as the optimization horizon is small enough.



\paragraph{Future research.} 
While \texttt{WFDroneBench} currently emphasizes optimization-based routing, reinforcement learning (RL)~\cite{electronics13132568} presents a promising avenue for adaptive decision-making in complex, uncertain environments. RL methods could enable end-to-end policies that respond directly to environmental feedback and stochastic fire behavior, potentially reducing dependence on precise real-time risk inputs. Other directions for future exploration include incorporating high-resolution, weather-conditioned dynamic risk maps, designing routing strategies that minimize expected detection time, and investigating heterogeneous fleet deployments. Additional work may also involve integrating dynamic predictions with real-time decision-making and expanding the benchmark suite with datasets derived from historical wildfire propagation events. Additionally, future work could extend the framework’s relevance to policy makers and mechanical engineers by enabling deeper analyses of resource-allocation decisions and by quantifying how drone hardware characteristics -- such as battery capacity and communication range -- shape achievable detection performance, as illustrated in  Figure \ref{fig:optimizationframework}.

Our findings highlight two important open challenges: (i) achieving rapid, reliable detection of small fires, and (ii) improving wildfire risk-map prediction, where the gap between ground-truth ignition patterns and available risk maps underscores a substantial opportunity for machine learning innovation.

\begin{ack}
    The second author of this paper is funded by the Netherlands Organisation for Scientific Research (NWO) under the Rubicon grant 019.233SG.010.
\end{ack}

\section*{Data Availability}
All benchmark datasets associated with this work are publicly available through the WFDroneBench dataset repository on Hugging Face:
\url{https://huggingface.co/datasets/MasterYoda293/DroneBench}.

\section*{Code Availability}
The code for generating the benchmark instances, running the sensor placement and drone routing algorithms, and reproducing the computational experiments is publicly available in the corresponding GitHub repository:
\url{https://github.com/RomainPuech/wildfire_drone_routing/tree/benchmark-release}.



{\raggedright
\bibliography{references}

\begin{thebibliography}{10}

\bibitem{f12040453}
B.~F. Azevedo, T.~Brito, J.~Lima, and A.~I. Pereira.
\newblock Optimum sensors allocation for a forest fires monitoring system.
\newblock {\em Forests}, 12(4), 2021.

\bibitem{WildfireMonitoringBailonruiz22}
R.~Bailon-Ruiz, A.~Bit-Monnot, and S.~Lacroix.
\newblock Real-time wildfire monitoring with a fleet of uavs.
\newblock {\em Robotics and Autonomous Systems}, 152:104071, 2022.

\bibitem{BAO2015100}
S.~Bao, N.~Xiao, Z.~Lai, H.~Zhang, and C.~Kim.
\newblock Optimizing watchtower locations for forest fire monitoring using
  location models.
\newblock {\em Fire Safety Journal}, 71:100--109, 2015.

\bibitem{WildfireDetectionBouguettaya}
A.~Bouguettaya, H.~Zarzour, A.~M. Taberkit, and A.~Kechida.
\newblock A review on early wildfire detection from unmanned aerial vehicles
  using deep learning-based computer vision algorithms.
\newblock {\em Signal Processing}, 190:108309, 2022.

\bibitem{10756651}
C.~C. Chan, S.~A. Alvi, X.~Zhou, S.~Durrani, N.~Wilson, and M.~Yebra.
\newblock A survey on iot ground sensing systems for early wildfire detection:
  Technologies, challenges, and opportunities.
\newblock {\em IEEE Access}, 12, 2024.

\bibitem{MCLPChurchRevelle}
R.~Church and C.~ReVelle.
\newblock The maximal covering location problem.
\newblock {\em Papers of the Regional Science Association}, 32(1):101--118,
  1974.

\bibitem{curtis2018classifying}
P.~G. Curtis, C.~M. Slay, N.~L. Harris, A.~Tyukavina, and M.~C. Hansen.
\newblock Classifying drivers of global forest loss.
\newblock {\em Science}, 361(6407):1108--1111, 2018.

\bibitem{DroneroutingDeLaFuente}
R.~De~la Fuente, M.~Aguayo, and C.~Contreras~Bolton.
\newblock An optimization-based approach for an integrated forest fire
  monitoring system with multiple technologies and surveillance drones.
\newblock {\em European Journal of Operational Research}, 313:435--451, March
  2024.

\bibitem{rs16214019}
Dipraj Debnath, Fernando Vanegas, Juan Sandino, Ahmad~Faizul Hawary, and Felipe
  Gonzalez.
\newblock A review of uav path-planning algorithms and obstacle avoidance
  methods for remote sensing applications.
\newblock {\em Remote Sensing}, 16(21), 2024.

\bibitem{electronics13132568}
Kubilay Demir, Vedat Tumen, Selahattin Kosunalp, and Teodor Iliev.
\newblock A deep reinforcement learning algorithm for trajectory planning of
  swarm uav fulfilling wildfire reconnaissance.
\newblock {\em Electronics}, 13(13), 2024.

\bibitem{diGiuseppe2025california}
F.~Di~Giuseppe, T.~Hewson, J.~McNorton, M.~Parrington, J.~Keune, and
  S.~El~Garroussi.
\newblock 2025 california wildfires: insights from ecmwf forecasts.
\newblock {\em ECMWF Science Blog}, February 2025.
\newblock Accessed May 5, 2025.

\bibitem{dillon2023whp}
G.~K. Dillon.
\newblock Wildfire hazard potential for the united states (270-m), version
  2023. 4th edition.
\newblock Forest Service Research Data Archive, Fort Collins, CO, 2023.
\newblock Updated 17 July 2024.

\bibitem{WildfireDetectionDLGhali22}
R.~Ghali, M.~A. Akhloufi, and W.~S. Mseddi.
\newblock Deep learning and transformer approaches for uav-based wildfire
  detection and segmentation.
\newblock {\em Sensors}, 22(5), 2022.

\bibitem{DroneroutingPostGhamry17}
K.~A. Ghamry, M.~A. Kamel, and Y.~Zhang.
\newblock Multiple uavs in forest fire fighting mission using particle swarm
  optimization.
\newblock In {\em 2017 International Conference on Unmanned Aircraft Systems
  (ICUAS)}, pages 1404--1409, 2017.

\bibitem{Gurobi}
{Gurobi Optimization}.
\newblock Gurobi optimizer reference manual 8.1.1, 2019.

\bibitem{hu2025coexposure}
J.~K. Hu, A.~Tri{\v{s}}ovi{\'c}, A.~Bakshi, D.~Braun, F.~Dominici, and J.~A.
  Casey.
\newblock Coexposure to extreme heat, wildfire burn zones, and wildfire smoke
  in the western us from 2006 to 2020.
\newblock {\em Science Advances}, 11(18):eadq6453, 2025.

\bibitem{DroneroutingJemmali}
M.~Jemmali, B.~M.~L. Kayed, W.~Boulila, H.~Amdouni, and M.~T. Alharbi.
\newblock Optimizing forest fire prevention: Intelligent scheduling algorithms
  for drone-based surveillance system.
\newblock {\em Procedia Computer Science}, 225:1562--1571, 2023.
\newblock 27th International Conference on Knowledge Based and Intelligent
  Information and Engineering Systems (KES 2023).

\bibitem{WildfireDetectionDLKim}
S.‐Y. Kim and A.~Muminov.
\newblock Forest fire smoke detection based on deep learning approaches and
  unmanned aerial vehicle images.
\newblock {\em Sensors}, 23(12), 2023.

\bibitem{WildfireDetectionDLKumar24}
A.~Kumar, A.~Perrusquía, S.~Al-Rubaye, and W.~Guo.
\newblock Wildfire and smoke early detection for drone applications: A
  light-weight deep learning approach.
\newblock {\em Engineering Applications of Artificial Intelligence},
  136:108977, 2024.

\bibitem{NEURIPS2024_02e978a2}
Y.~Li, K.~Li, G.~Li, Z.~Wang, C.~Ji, L.~Wang, D.~Zuo, Q.~Guo, F.~Zhang,
  M.~Wang, and D.~Lin.
\newblock Sim2real-fire: A multi-modal simulation dataset for forecast and
  backtracking of real-world forest fire.
\newblock {\em Advances in Neural Information Processing Systems},
  37:1428--1442, 2024.

\bibitem{mugin_ev350}
{Mugin UAV}.
\newblock {Mugin EV350 (KIT) Full Electric Carbon Fiber VTOL UAV Platform}.
\newblock
  \url{https://www.muginuav.com/product/mugin-ev350-carbon-fiber-full-electric-vtol-uav-platform/},
  2025.
\newblock Accessed: 2025-05-15.

\bibitem{drones9040263}
Mamunur Rahman, Nurul~I. Sarkar, and Raymond Lutui.
\newblock A survey on multi-uav path planning: Classification, algorithms, open
  research problems, and future directions.
\newblock {\em Drones}, 9(4), 2025.

\bibitem{s23031463}
Asif~Mahmud Raivi, S.~M.~Asiful Huda, Muhammad~Morshed Alam, and Sangman Moh.
\newblock Drone routing for drone-based delivery systems: A review of
  trajectory planning, charging, and security.
\newblock {\em Sensors}, 23(3), 2023.

\bibitem{RAMADAN2024101248}
M.N.A. Ramadan, T.~Basmaji, T.~Gad, H.~Hamdan, B.T. Akgün, M.A.H. Ali,
  M.~Alkhedher, and M.~Ghazal.
\newblock Towards early forest fire detection and prevention using ai-powered
  drones and the iot.
\newblock {\em Internet of Things}, 27:101248, 2024.

\bibitem{roque2025identifying}
N.~A. Roque, H.~Andrews, and A.~R. Santos-Lozada.
\newblock Identifying air quality monitoring deserts in the united states.
\newblock {\em Proceedings of the National Academy of Sciences},
  122(17):e2425310122, 2025.

\bibitem{DroneroutingSantin}
R.~Santin, L.~Assis, A.~Vivas, and L.~C.~A. Pimenta.
\newblock Matheuristics for multi-uav routing and recharge station location for
  complete area coverage.
\newblock {\em Sensors}, 21(5), 2021.

\bibitem{9068226}
J.~Shi, W.~Wang, Y.~Gao, and N.~Yu.
\newblock Optimal placement and intelligent smoke detection algorithm for
  wildfire-monitoring cameras.
\newblock {\em IEEE Access}, 8:72326--72339, 2020.

\bibitem{short2021spatial}
K.~C. Short.
\newblock Spatial wildfire occurrence data for the united states, 1992–2020
  [fpa\_fod\_20221014], 2022.
\newblock U.S. Department of Agriculture, Forest Service.

\bibitem{thomas2017costs}
D.~Thomas, D.~Butry, S.~Gilbert, D.~Webb, and J.~Fung.
\newblock The costs and losses of wildfires: A literature survey (nist special
  publication 1215).
\newblock {\em National Institute of Standards and Technology. U.S. Department
  of Commerce}, 2017.

\bibitem{tyukavina2022global}
A.~Tyukavina, P.~Potapov, M.~C. Hansen, A.~H. Pickens, S.~V. Stehman,
  S.~Turubanova, D.~Parker, V.~Zalles, A.~Lima, I.~Kommareddy, et~al.
\newblock Global trends of forest loss due to fire from 2001 to 2019.
\newblock {\em Frontiers in Remote Sensing}, 3:825190, 2022.

\bibitem{wang2021economic}
D.~Wang, D.~Guan, S.~Zhu, M.~M. Kinnon, G.~Geng, Q.~Zhang, H.~Zheng, T.~Lei,
  S.~Shao, P.~Gong, et~al.
\newblock Economic footprint of california wildfires in 2018.
\newblock {\em Nature Sustainability}, 4(3):252--260, 2021.

\bibitem{fire7070254}
Y.~Wang, F.~Gao, and M.~Li.
\newblock Probabilistic path planning for uavs in forest fire monitoring:
  Enhancing patrol efficiency through risk assessment.
\newblock {\em Fire}, 7(7), 2024.

\bibitem{DroneroutingXu}
Y.~Xu, J.~Li, and F.~Zhang.
\newblock A uav-based forest fire patrol path planning strategy.
\newblock {\em Forests}, 13(11), 2022.

\bibitem{WildfireriskXu2024}
Z.~Xu, J.~Li, S.~Cheng, X.~Rui, Y.~Zhao, H.~He, and L.~Xu.
\newblock Wildfire risk prediction: A survey of recent advances using deep
  learning techniques.
\newblock arXiv preprint arXiv:2405.01607, 2024.

\bibitem{DroneroutingYan}
X.~Yan and R.~Chen.
\newblock Application strategy of unmanned aerial vehicle swarms in forest fire
  detection based on the fusion of particle swarm optimization and artificial
  bee colony algorithm.
\newblock {\em Applied Sciences}, 14(11), 2024.

\bibitem{WildfireDetectionDLYandouzi22}
M.~Yandouzi, M.~Grari, I.~Idrissi, M.~Boukabous, O.~Moussaoui, M.~Azizi,
  K.~Ghoumid, and A.~K. Elmiad.
\newblock Forest fires detection using deep transfer learning.
\newblock {\em Forest}, 13(8):1, 2022.

\bibitem{WildfireDetectionDLYang}
H.~Yang, J.~Wang, and J.~Wang.
\newblock Efficient detection of forest fire smoke in uav aerial imagery based
  on an improved yolov5 model and transfer learning.
\newblock {\em Remote Sensing}, 15(23), 2023.

\end{thebibliography}
}

\bibliographystyle{plain}






\appendix

\section{Summary of sets, parameters, and decision variables}\label{appendix:sets}

To ensure clarity in the mathematical formulations presented in the appendices, Table~\ref{tab:overview} summarizes all sets, parameters, and decision variables referenced in the appendices.
\small
\begin{longtable}{lll}
\caption{Overview of sets, parameters, and decision variables used in the supplementary materials.} \label{tab:overview} \\
\toprule
\textbf{Sets} & \textbf{Description} & \textbf{Definition} \\
\midrule
\endfirsthead

\multicolumn{3}{c}%
{{\bfseries \tablename\ \thetable{} -- continued from previous page}} \\
\toprule
\textbf{Parameters} & \textbf{Description}  \\
\midrule
\endhead

\midrule \multicolumn{3}{r}{{Continued on next page}} \\
\endfoot

\bottomrule
\endlastfoot

$\set I$ & Set of all possible grid points & \\
$\set I_g$ & Feasible grid points for ground stations & $\subseteq \set I$ \\
$\set I_c$ & Feasible grid points for charging stations & $\subseteq \set I$ \\
$\set I_d$ & Feasible grid points for drones & $\subseteq \set I$ \\
$\set A_d(i)$ & Neighboring grid points of $i$ within \\
& distance $d$ (including $i$) & \\
$\set T$ & Set of all time periods & $\{1, \ldots, T\}$ \\
$\set H$ & Time periods in optimization time horizon & $\{1, \ldots, H\}$ \\
$\set S$ & Set of drones & $\{1, \ldots, S\}$ \\
$\set C$ & Set of charging stations & Cardinality $C$ \\
$\set G$ & Set of ground sensors & Cardinality $G$ \\
$\set N_{c}$ & Nearby charging station pairs & $\{(i,j) \in \set I_c^2 \mid i \neq j, \ \norm{i-j} \leq \rho_c\}$ \\
$\set N_{g}$ & Nearby ground sensor pairs & $\{(i,j) \in \set I_g^2 \mid i \neq j, \ \norm{i-j} \leq \rho_g\}$ \\
$\set N_{cg}$ & Charging–ground sensor pairs within & $\{(i,j) \in \set I_c \times I_g \mid i \neq j,$ \\
& exclusion radius & $\quad \norm{i-j} \leq \rho_{cg}\}$\\
$\set W$ & Set of historical wildfires \\
\midrule
\textbf{Parameters} & \multicolumn{2}{l}{\textbf{Description}} \\
\midrule
$\Delta_i$ & \multicolumn{2}{l}{Minimum distance from grid point $i$ to nearest charging station} \\
$B_{\text{max}}$ & \multicolumn{2}{l}{Maximum number of time steps a drone can fly without recharging} \\
$r_{i}^s$ & \multicolumn{2}{l}{Static wildfire risk at grid point $i$} \\
$r_{it}^d$ & \multicolumn{2}{l}{Dynamic wildfire risk at grid point $i$ at time $t$} \\
$f_{it}^w$ & \multicolumn{2}{l}{Indicator of if grid point $i$ is burnt at time $t$ in wildfire scenario $w$.} \\
$n_g$ & \multicolumn{2}{l}{Maximum number of ground sensors that can be placed} \\
$n_c$ & \multicolumn{2}{l}{Maximum number of charging stations that can be placed} \\
$n_d$ & \multicolumn{2}{l}{Maximum number of drones assignable to each charging station} \\
$\tau$ & \multicolumn{2}{l}{Coverage memory window: revisits provide no reward during this time} \\
$\rho_c$ & \multicolumn{2}{l}{Exclusion radius: minimum spacing between charging stations} \\
$\rho_g$ & \multicolumn{2}{l}{Exclusion radius: minimum spacing between ground sensors} \\
$\rho_{cg}$ & \multicolumn{2}{l}{Exclusion radius: minimum spacing between charging stations and ground sensors} \\
$\lambda_{\text{disp}}$ & \multicolumn{2}{l}{Weight of drone dispersion in the objective function} \\
$t_{\text{det}}^w$ & \multicolumn{2}{l}{Time of detection of wildfire $w \in \set W$} \\
$t_{\text{ign}}^w$ & \multicolumn{2}{l}{Time of ignition of wildfire $w \in \set W$} \\
$t_{\text{det}}^s$ & \multicolumn{2}{l}{Time at which drone $s$ detects a wildfire; equals $T$ if no fire is detected by drone $s$} \\
$t_{\text{exe}}$ & Execution time of the drone routing algorithm \\ 
$\bm p^c$ & Grid coordinates in $\mathbb{Z}^2$ of charging station $c$ \\ 
$T_{\text{reev}}$ & Reevaluation step for rolling optimization \\ 
    $D_i$ & \multicolumn{2}{l}{\blue{Distance from grid point i to the nearest charging station}} \\ \midrule
\textbf{Variables} & \multicolumn{2}{l}{\textbf{Description}} \\ \midrule
    $x_i^g$ & \multicolumn{2}{l}{\text{Binary variable denoting if a ground sensor is placed at grid point i}} \\
    $x_i^c$ & \multicolumn{2}{l}{\text{Binary variable denoting if a charging station is placed at grid point i}} \\
    $a_{its}$ & \multicolumn{2}{l}{\text{Binary variable denoting if drone $s$ flies at grid point $i$ at time $t$}} \\
    $c_{its}$ & \multicolumn{2}{l}{\text{Binary variable denoting if drone $s$ charges at grid point $i$ at time $t$}} \\
    $\theta_{it}$ & \multicolumn{2}{l}{\text{Binary variable denoting if grid point $i$ is covered by a drone at time $t$}} \\
    $b_{ts}$ & \multicolumn{2}{l}{\text{Battery of drone $s$ at time $t$ defined as the number of time steps drone $s$ can operate}} \\
    & \multicolumn{2}{l}{\text{without recharging}} \\
    $\bm p^{s,t}$ & \multicolumn{2}{l}{\text{Grid coordinates in $\mathbb{Z}^2$ of drone $s$ at time $t$}} \\
\end{longtable}
\normalsize

\section{Definitions of evaluation metrics}\label{appendix:evaluation metrics}

In Table \ref{tab:evaluation metrics} we give the mathematical formulations of the evaluation metrics as described in Section \ref{subsec:evaluationmetrics} of the main paper. 
\begin{table}[!htb]
\centering
\resizebox{0.95\textwidth}{!}{%
\begin{tabular}{lll} \toprule
\textbf{Evaluation metric}             & \textbf{Formula}                  \\ \midrule
Mean detection time (MDT) & $\frac{1}{\sum_{w \in \mathcal{W}} \mathds{1}\{t_{\text{det}}^w \leq T \}} \sum_{w \in \mathcal{W}} \mathds{1}\{t_{\text{det}}^w \leq T \} (t_{\text{det}}^w - t_{\text{ign}}^w)$ \\ \addlinespace
Mean fire size at detection (MSD) & $\frac{1}{|\mathcal{I}|}\sum_{i \in \mathcal{I}}f_{it}^w, \ \text{with } t = t_{\text{exe}}^w$ \\ \addlinespace
\% \text{ fires detected} (FDR) &  $\frac{1}{|\set W|} \sum_{w \in \set W} \mathds{1}\{t_{\text{det}}^w \leq T \} \cdot 100 \%$ \\ \addlinespace
\text{Proportion of detection by device D} (DS-D) & $\displaystyle \frac{\sum_{w \in \set W} \mathds{1}\{t_{\text{det}}^w \leq T \} \mathds{1}\{\text{fire detected by device D}\}}{\sum_{w \in \set W} \mathds{1}\{t_{\text{det}}^w \leq T \}}$ \\ \addlinespace
Mean percentage of map explored (MCE) & 
    $\frac{1}{|\set I|} \sum_{i \in \set I} \mathds{1}\left(\sum_{t \in \set T} \theta_{it} \geq 1 \vee x_i^g = 1 \vee x_i^c = 1 \right) $ \\ \addlinespace
Mean distance traveled by drones (MDist) & $\frac{1}{S} \sum_{s \in \set S} \sum_{t=1}^{t_{\text{det}}^s} d\left(\bm{p^{s,t}},\bm{p^{s,t+1}}\right)$ \\ \addlinespace
Mean distance to nearest charging station (DCSD) & $\frac{1}{ST} \sum_{s \in \set S} \sum_{t \in \set T} \min_{c \in \set C} d\left(\bm{p^{s,t}},\bm{p^c} \right)$ \\ \addlinespace
Mean drone routing execution time (RTE) & $\frac{1}{\# \text{scenarios}} \sum_{i=1}^{\# \text{scenarios}} t_{\text{exe}}$  \\
\bottomrule \addlinespace
\end{tabular}}
\caption{Mathematical formulations of the evaluation metrics.}
\label{tab:evaluation metrics}
\end{table}

\section{Decomposition algorithm}\label{appendix:decomposition alg}
To enhance scalability and reduce computational complexity, we partition the full set of charging stations into spatial operational zones, where each operational zone includes stations that are mutually reachable within a drone’s battery range. Specifically, two stations are placed in the same operational zone if the Euclidean distance between them is less than or equal to the drone's battery range, ensuring that any drone can travel between stations in the operational zone. The decomposition algorithm is given in Algorithm \ref{alg:clustering_model}.
 \begin{algorithm}[h]
        \caption{Decomposition model \label{alg:clustering_model}} 
        \textbf{Input}: set of charging stations $\set C$. \\ 
        \textbf{Output}: operational zones $\set K_1, \ldots, \set K_m$, where $\set K_k \subseteq \set C$ for all $k \in [m]$, $\cap_{k \in [m]} \set K_k = \emptyset$.

      \begin{algorithmic}[1]
            \STATE Set $k = 1$ 
            \WHILE{$\set C \neq \emptyset$}
            \STATE Select a charging station $i \in \set C$
            \STATE Remove $i$ from $\set C$
            \STATE Set $\set K_k = \{i\}$
            \FOR{$j \in \set C \setminus \{i\}$}
            \IF{the Euclidean distance between $i$ and $j$ is less or equal than $\lambda$}
            \STATE add $j$ to operational zone $\set K_k$
            \STATE remove $j$ from $\set C$
            \ENDIF
            \ENDFOR
            \ENDWHILE
            \RETURN $\set K_1, \ldots, \set K_m$.
         \end{algorithmic}
    \end{algorithm}

\section{Parameter values used in the experiments}\label{appendix:par values}
In Table \ref{tab:parameter_values} we give an overview of the parameter values used in the experiments when using the max-coverage ground sensor placement baseline in combination with a drone routing strategy.
\begin{table}[H]
    \centering
        \caption{Parameter values used in the experiments.}
    \begin{tabular}{lc}
    \toprule
    \textbf{Parameter} & \textbf{Value} \\ \midrule
    $n_g$ & 8 \\
    $n_c$ & 2  \\
    $n_d$ & 2  \\
    $\rho_c$ & 10 \\
    $\rho_{cg}$ & 10 \\
    $\rho_g$ & 1 \\
    Drone speed (m/min) & 600  \\
    Coverage radius (m) & 300 \\
    Grid point size (m) & 30 \\
    Transmission range (m) & 50000 \\
    $T_{\text{reev}}$ & 5  \\
    $H$ & 10  \\
    $B_{\text{max}}$ (h) & 1  \\
    \bottomrule
    \end{tabular}
    \label{tab:parameter_values}
\end{table}

\section{Mathematical model formulations}\label{appendix:math model formulations}
In this section the mathematical formulations are given for all models. All distances are measured in the $L^{\infty}$ norm, unless specified otherwise.

\subsection{Ground sensors and charging stations placement model}\label{appendix:max coverage ground}
The simplest version of the ground sensors and charging stations max-coverage model only considers direct coverage of a device at its immediate location. \textit{We first present this model before extending it to the indirect and partial coverage model we used.} The simple model is given by:
\allowdisplaybreaks
\begin{subequations}\label{eq:maxcov_placement}
    \begin{align} 
        \max_{\bm x^g, \bm x^c} &\quad \displaystyle \sum_{i \in \set I_g} r_i^s x_i^g + \sum_{i \in \set I_c} r_i^s x_i^c \label{eq:maxcov_placement_a} \\
        {\rm s.t. } &\quad \displaystyle x_i^g + x_i^c \leq 1, &\quad \forall i \in \set I_g \cap \set I_c, \label{eq:maxcov_placement_b} \\
        &\quad \displaystyle \sum_{i \in \set I_g} x_i^g \leq n_g, \label{eq:maxcov_placement_c} \\
        &\quad \displaystyle \sum_{i \in \set I_c} x_i^c \leq n_c, \label{eq:maxcov_placement_d} \\
        &\quad \displaystyle x_i^g + x_i^c \leq 1, &\quad (i,j) \in \set N_g, \label{eq:maxcov_placement_e} \\
        &\quad \displaystyle x_i^g + x_j^c \leq 1, &\quad (i,j) \in \set N_c, \label{eq:maxcov_placement_f} \\
        &\quad \displaystyle x_j^g + x_i^c \leq 1, &\quad (i,j) \in \set N_{cg}.\label{eq:maxcov_placement_g} 
    \end{align}
\end{subequations}
Here, constraints \eqref{eq:maxcov_placement_b} ensure that there can either be a charging station or a ground sensor at grid point $i$, but not both. Constraints \eqref{eq:maxcov_placement_c} and \eqref{eq:maxcov_placement_d} are capacity constraints on the ground sensors and charging stations respectively. Finally, \eqref{eq:maxcov_placement_e}-\eqref{eq:maxcov_placement_g} are spatial exclusion constraints between two ground sensors, two charging stations, and a ground sensor and a charging station, respectively.

\textbf{Extensions.}
Observe that we have assumed that a ground sensor as well as a charging station only detects a wildfire at the grid point it occupies. This can be naturally extended to allow partial detection of wildfires in neighboring grid points by replacing \eqref{eq:maxcov_placement} by 
\begin{align*}
    \displaystyle \max_{\bm \phi, \bm x, \bm y} &\quad \displaystyle \sum_{i \in \set I_g \cup \set I_c} \phi_i \\
    {\rm s.t. } 
    &\quad \displaystyle \phi_i \leq \sum_{k \in \set I_g \cap \set A_d(k)} \gamma_{ik}^g x_k^g + \sum_{k \in \set I_c \cap \set A_{d'}(k)} \gamma_{ik}^c x_k^c, &\quad i \in \set I_g \cup \set I_c,\\
    &\quad \displaystyle 0 \leq \phi_i \leq 1, &\quad i \in \set I_g \cup \set I_c,
\end{align*}
where $\gamma_{ik}^G$ and $\gamma_{ik}^C$ denote the partial contribution from a ground sensor and charging station at grid point $k$ to grid point $i$ respectively, and $\phi_i^g$ denotes the total fractional coverage in grid point $i$ due to ground sensors and charging stations nearby.


Observe that our model readily allows for the integration of existing ground sensors alongside newly optimized placements by explicitly including them in the formulation with fixed variable values (e.g., setting $x_i^g = 1$ for grid points with pre-existing sensors) or alternatively excluding them from the decision variable domains.

\textbf{Gaussian Coverage}
Notice that the Gaussian coverage model we use and presented in section \ref{gaussian-coverage} is a special case of the extension presented above with the following parameters: 
\begin{itemize}
    \item $\gamma^g_{ik}$ is 1 for $i = k$ and 0 otherwise, as ground sensors only cover their immediate location.
    \item All $\gamma^c_{ik}$ are precomputed at the same time for a given $i$. We start with a grid of zeros with a unique 1 at cell $i$ and iteratively take convolutions with a $3 \times 3$ kernel of ones, simulating the diffusion process of a drone over neighboring cells. After dividing the resulting grid by the number of convolutions applied, the value of cell $k$ gives the expected number of visits of cell $k$ by a drone following a Brownian motion starting at $i$. Capping this value by 1 gives $\gamma^c_{ik}$.
\end{itemize}
 Expressing the Gaussian coverage as a special case of the extended max coverage model allows us to solve it efficiently using linear programming.

\subsection{Max coverage drone routing model}\label{appendix:max coverage drones}
The drone routing optimization model with the objective to maximize coverage is given by:
\allowdisplaybreaks
\begin{subequations} \label{eq:maxcov}
    \begin{align}
        \max_{\bm a, \bm b, \bm c, \bm \theta}  &\quad \displaystyle \sum_{i \in \set I_d} \left\{r_{i1}^d \theta_{i1} + \sum_{t \in \set H \setminus \{1\}} r_{it}^d (\theta_{it} - \theta_{i,t-1}) \right\} \label{eq:maxcov_a}
        \\ {\rm s.t. } 
        &\quad \displaystyle \sum_{i \in \set I_d} a_{its} + \sum_{i \in \set C} c_{its} = 1, & t \in \set H \setminus \{1\}, s \in \set S, \label{eq:maxcov_b}\\
        &\quad \displaystyle \sum_{s \in \set S} c_{its} \leq n_d, & i \in \set C, t \in \set H \setminus \{1\}, \label{eq:maxcov_c}\\
        &\quad \displaystyle c_{j,t+1,s} + a_{j,t+1,s} \leq \sum_{i \in \set I_d \cap \set A_1(j)} a_{its} + c_{jts}, &  j \in \set C, t \in \set H \setminus \{H\}, s \in \set S, \label{eq:maxcov_d}\\
        &\quad \displaystyle a_{j,t+1,s} \leq \sum_{i \in \set I_d \cap \set A_1(j)} a_{its}, &  j \in \set I_d \setminus \set C, t \in \set H \setminus \{H\}, s \in \set S, \label{eq:maxcov_e}\\
        &\quad \displaystyle 0 \leq b_{ts} \leq B_{\text{max}}, & t \in \set H, s \in \set S, \label{eq:maxcov_f}\\
        &\quad \displaystyle b_{ts} \geq B_{\text{max}} \sum_{i \in \set C} c_{its}, & t \in \set H, s \in \set S, \label{eq:maxcov_g}\\
        &\quad \displaystyle b_{t+1,s} \leq b_{ts} - 1 + (B_{\text{max}} + 1) \sum_{j \in \set C} c_{j,t+1,s} & t \in \set H \setminus \{H\}, s \in \set S, \label{eq:maxcov_h}\\
        &\quad \displaystyle b_{Ts} \geq D_i \cdot a_{iTs}, & i \in \set I_d, s \in \set S, \label{eq:maxcov_i}\\
        &\quad \displaystyle \theta_{it} = 0, & i \in \set I_d \cap (\set G \cup \set C), t \in \set H, \label{eq:maxcov_l} \\
        &\quad \displaystyle \theta_{it} \geq a_{its}, &\quad i \in \set I_d, t \in \set T, s \in \set S, \label{eq:maxcov_m} \\
        &\quad \displaystyle \theta_{1t} \leq \sum_{s \in \set S} a_{1ts}, &\quad t \in \set H, \label{eq:maxcov_n} \\
        &\quad \displaystyle \theta_{it} \leq \sum_{s \in \set S} a_{its} + \theta_{i,t-1}, &\quad i \in \set I_d, t \in \set H \setminus \{1\}, \label{eq:maxcov_o} \\
        &\quad \displaystyle \theta_{it} \geq \theta_{i,t-1}, &\quad i \in \set I_d, t \in \set H \setminus \{1\}. \label{eq:maxcov_p}
    \end{align}
\end{subequations}
Constraints \eqref{eq:maxcov_b} ensure that each drone at time $t$ is either charging or flying, but not both. 
Constraints \eqref{eq:maxcov_c} ensure that there can be at most $n_d$ drones at a time at a charging station.
Constraints \eqref{eq:maxcov_d} - \eqref{eq:maxcov_e} ensure that a drone can only fly or charge at location $j$ at time $t+1$ if it was charging already in the same location or the drone was in a neighboring location at time $t$.
\eqref{eq:maxcov_f}  - \eqref{eq:maxcov_h} are battery constraints that ensure that when a drone is charging, the battery updates to $B_{\text{max}}$, when a drone is flying, the battery decreases with 1 every time step, and the battery does not drop below $0$. Constraints \eqref{eq:maxcov_i} enforce a ``no-suicide" condition for drones at the end of the planning horizon, that is, a drone must have enough battery remaining at the final time step to reach the nearest charging station from its final position. Constraints \eqref{eq:maxcov_l}-\eqref{eq:maxcov_p} are coverage constraints. 
Constraints \eqref{eq:maxcov_l} make sure that there is no need to consider drone-based coverage at points where ground sensors or charging stations are already installed. Constraints \eqref{eq:maxcov_m} ensures that grid point $i$ must be covered at time $t$ when there is a drone present at grid point $i$ at time $t$. Constraints \eqref{eq:maxcov_n} and \eqref{eq:maxcov_o} ensure that a grid point is only covered at time $t$ if at least one drone visits it at time $t$ or it was already covered at time $t-1$. Constraints \eqref{eq:maxcov_p} ensure that once a grid point is detected, it remains detected in all future time steps.

\begin{figure}[htbp]
    \centering
\includegraphics[width=0.95\textwidth]{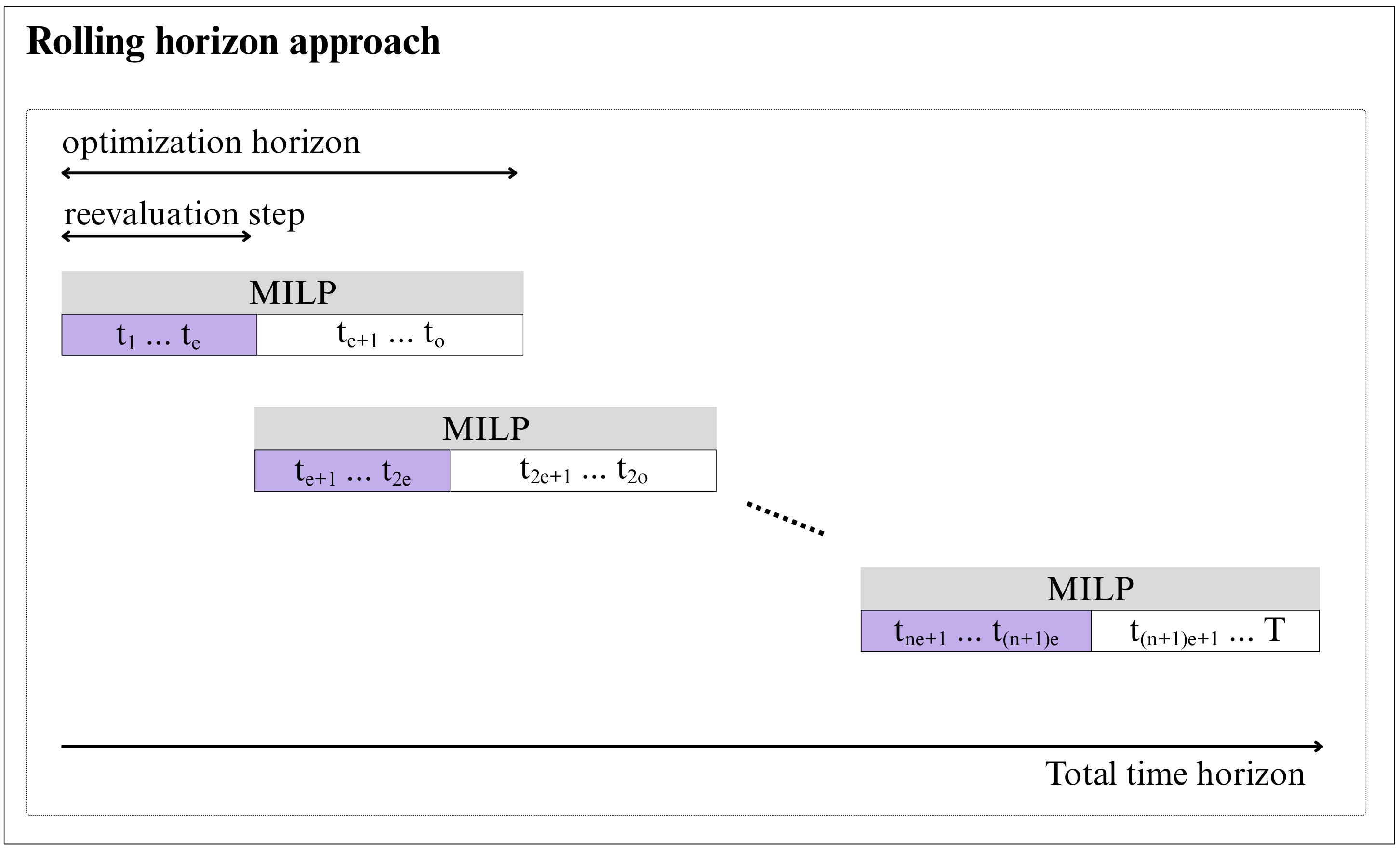}
    \caption{The rolling horizon optimization framework. The problem is solved over an optimization horizon of $o$ time steps using Mixed Integer Linear Programming (MILP). Even though solutions are generated for $o$ time steps, we only use the $e$ first ones, with $e$ the reevaluation step parameter. After each reevaluation step of $e$ time steps, the horizon is shifted forward, and the problem is solved again. This process continues until the end of the total time horizon is reached, ensuring adaptability to new information and improved solution quality over time.}
    \label{fig:rollinghorizon}
\end{figure}

The drone routing is performed using a rolling-horizon optimization framework (see Figure \ref{fig:rollinghorizon}). At each decision step, we solve the above MILP over a limited optimization horizon to account for the future implications of routing decisions while maintaining tractability. Only the immediate $T_{\text{reev}}$ steps (depicted in purple in Figure \ref{fig:rollinghorizon}) are executed before risk data is updated and the optimization is re-solved. This enables the system to adapt dynamically to evolving wildfire conditions based on drone observations (e.g., risk decay following surveillance) or external updates such as weather changes. 

In addition, we initialize the drone routing model by adding the following constraints to the optimization problem formulation
\begin{align}
    &\quad \displaystyle \sum_{i \in \set C} a_{i1s} + \sum_{i \in \set C} c_{i1s} = 1, &  s \in \set S, \label{eq:maxcov_b_init}\\
    &\quad \displaystyle \sum_{s \in \set S} (a_{i1s} + c_{i1s}) \leq n_d, & i \in \set C, \label{eq:maxcov_d_init}\\
    &\quad \displaystyle b_{1s} = B_{\text{max}}, & s \in S, \label{eq:maxcov_h_init}
\end{align}
ensuring that all drones start flying from a charging station (constraints \eqref{eq:maxcov_b_init}), there are at most $n_d$ drones either flying or charging at a charging station at $t = 1$ (constraints \eqref{eq:maxcov_d_init}), and each drone starts with full capacity (constraints \eqref{eq:maxcov_h_init}). 

\textbf{Extensions.}
Observe that we have assumed that a drone only detects a wildfire at the grid point it occupies. This can be naturally extended to allow partial detection of wildfires in neighboring grid points by replacing constraint \eqref{eq:maxcov_m} by 
\begin{align*}
            &\quad \displaystyle \theta_{it} \leq \sum_{s \in S} \sum_{k \in \set I_d \cap \set A_d(i)}\gamma_{ik}^d a_{kts}, &\quad i \in \set I_d, t \in \set T,
\end{align*}
where $\gamma_{ik}^d$ denotes the partial contribution from a drone at grid point $k$ to grid point $i$. In this case, the variable $\theta_{it}$ 
may take fractional values and is interpreted as the degree of coverage at location $i$ and time $t$.

To encourage drones to spatially disperse we can include a separation reward term in the objective, defined as the sum of pairwise $L_{\infty}$ distances between drones over time, i.e., we consider \eqref{eq:maxcov} in which we replace the objective \eqref{eq:maxcov_a} by
\begin{align} \label{eq:maxcov_a_reg}
    \max \left\{ \sum_{i \in \set I_d} \sum_{t \in \set T} r_{it} \phi_{it} + \lambda_{\text{disp}} \sum_{t \in \set T} \sum_{s_1 \in \set S} \sum_{s_2 < s_1, s_2 \in \set S} \max \left\{\big|p^{s_1,t}_1 - p^{s_2,t}_1\big|, \big|p^{s_1,t}_2 - p^{s_2,t}_2\big| \right\} \right\}. 
    \end{align}
Note that we can easily linearize the objective \eqref{eq:maxcov_a_reg} using auxiliary variables and big-M constraints.

\section{Results of remaining evaluation metrics}\label{appendix:evaluation metrics results}
The results across strategies and layouts for each secondary metric are displayed in Table~\ref{additional_metrics_bm} and Table~\ref{additional_metrics_bp}.

\begin{table}[htbp]
\centering
\begin{minipage}{0.48\textwidth}
\centering
\rotatebox{90}{%
\begin{minipage}{.85\textheight} 
\centering
\resizebox{\textwidth}{!}{%
\begin{tabular}{llllllll}
\toprule
 &  & execution time & fire size cells & fire percentage & map explored & total distance & drone detection pct \\
\midrule
sensor strategy & drone strategy &  &  &  &  &  &  \\
\midrule
\midrule
\multirow[t]{4}{*}{Random} & MaxCov & 20.315 (9.243) & 317.865 (500.740) & 0.115 (0.162) & 0.089 (0.058) & 10811.368 (8586.738) & 99.7\% (1.0\%) \\
\midrule
 & UniCov & 13.777 (1.405) & 394.954 (659.147) & 0.137 (0.220) & 0.133 (0.051) & 12539.380 (8558.823) & 100.0\% (0.0\%) \\
\midrule
 & Brownian & 0.668 (0.150) & 1891.451 (2025.248) & 0.818 (0.699) & 0.346 (0.224) & 40051.673 (19738.856) & 80.6\% (37.9\%) \\
\midrule
\cline{1-8}
\multirow[t]{4}{*}{GaussianCov} & MaxCov & 19.748 (9.022) & 198.689 (239.398) & 0.076 (0.076) & 0.106 (0.058) & 11170.125 (6990.986) & 94.0\% (14.1\%) \\
\midrule
 & UniCov & 12.263 (1.566) & 295.289 (387.528) & 0.107 (0.127) & 0.120 (0.063) & 14744.375 (13057.346) & 95.2\% (10.4\%) \\
\midrule
 & Brownian & 0.643 (0.133) & 2011.373 (2886.966) & 0.835 (1.081) & 0.346 (0.213) & 39164.918 (20949.228) & 71.5\% (40.2\%) \\
\midrule
\cline{1-8}
\bottomrule
\end{tabular}}
\caption{Results with baseline sensor strategies}\label{additional_metrics_bm}
\end{minipage}}
\end{minipage}
\hfill
\begin{minipage}{0.48\textwidth}
\centering
\rotatebox{90}{%
\begin{minipage}{.85\textheight}
\centering
\resizebox{\textwidth}{!}{%
\begin{tabular}{llllllll}
\toprule
 &  & execution time & fire size cells & fire percentage & map explored & total distance & drone detection pct \\
\midrule
sensor strategy & drone strategy &  &  &  &  &  &  \\
\midrule
\midrule
\multirow[t]{4}{*}{ Random} & MaxCov & 0.000 (0.000) & 641.130 (1152.963) & 0.225 (0.396) & 0.113 (0.060) & 14930.118 (12428.318) & 99.4\% (1.1\%) \\
\midrule
 & UniCov & 0.000 (0.000) & 577.594 (1180.993) & 0.204 (0.407) & 0.132 (0.057) & 14472.535 (11741.451) & 99.5\% (1.0\%) \\
\midrule
 & Brownian & 0.618 (0.120) & 1094.617 (1523.182) & 0.367 (0.489) & 0.130 (0.070) & 18401.000 (16735.612) & 97.4\% (7.7\%) \\
\midrule
\cline{1-8}
\multirow[t]{4}{*}{GaussianCov} & MaxCov & 0.000 (0.000) & 256.074 (301.431) & 0.098 (0.098) & 0.123 (0.053) & 12036.871 (7137.082) & 99.7\% (0.9\%) \\
\midrule
 & UniCov & 0.000 (0.000) & 237.247 (299.693) & 0.086 (0.098) & 0.116 (0.031) & 10157.665 (6907.177) & 100.0\% (0.1\%) \\
\midrule
 & Brownian & 0.627 (0.126) & 1029.781 (1240.544) & 0.372 (0.400) & 0.206 (0.130) & 24515.343 (15633.415) & 96.2\% (9.6\%) \\
\midrule
\cline{1-8}
\bottomrule
\end{tabular}}
\caption{Results with Gaussian sensor strategies (Table 2)}\label{additional_metrics_bp}
\end{minipage}}
\end{minipage}
\end{table}

\section{Licensing datasets \& assets used}

We use several publicly available datasets in this study and ensure that each is properly credited and used in accordance with its license and usage terms. For the Sim2Real-Fire dataset, we use simulation data, which is released under the Apache-2.0 license as noted in the original publication and hosted repository (available at \url{https://github.com/TJU-IDVLab/Sim2Real-Fire}). The Burn Probability (BP) for the United States (270-m), version 2023 (4th Edition) dataset, released by the USDA Forest Service, is publicly accessible for research and other uses under a Creative Commons Attribution (CC BY) license. It is available at \url{https://www.fs.usda.gov/rds/archive/}. We also incorporate the Fire Program Analysis Fire-Occurrence Database (FPA FOD), maintained by the U.S. Forest Service, which is available for public research under open data access policy (CC0 1.0 License). All cited datasets are accompanied by their corresponding original publications in our bibliography. 

A summary of the asset licenses and their corresponding source web locations is provided in Table \ref{tab:assetlicenses}. 

{
\small
\begin{longtable}{p{4cm}p{6cm}p{2.5cm}}
\toprule
\textbf{Dataset} & \textbf{URL or DOI}                               & \textbf{License}   \\
\endfirsthead
\endhead
\midrule
National USFS Fire Occurrence Point (Feature Layer) &
  \url{https://data-usfs.hub.arcgis.com/datasets/6059c1a4dca749d393e33ee5f8a0cbaf_9/about} &
  CC0 1.0 \\ \midrule
Sim2Real-Fire    & \url{https://github.com/TJU-IDVLab/Sim2Real-Fire} & Apache-2.0 \\ \midrule
Wildfire Hazard Potential for the United States (270-m), version 2023 (4th Edition) &
  \url{10.2737/RDS-2015-0047-4} &
  Creative \newline Commons \newline CC-BY
  \\ \bottomrule
  \caption{Overview of the used asset licenses and their corresponding source web locations. \label{tab:assetlicenses}} \\
\end{longtable}
}

\section{Instructions to download dataset and code}\label{appendix: instructions to download dataset and code}
Our dataset and code are available in the supplementary material. Please refer to the README file for how to run the code and reproduce the results of our experiments.


\section{Dataset plots}\label{appendix:dataset_plots}
\begin{figure}[htbp]
    \includegraphics[width=\textwidth]{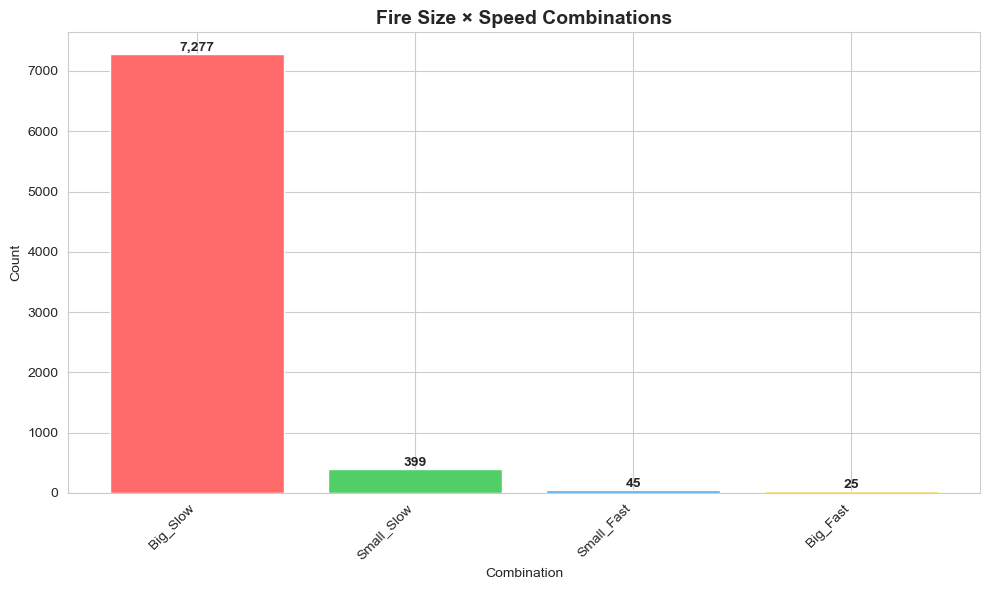}\label{fig:size_and_speed}
    
    \caption{\rev{Distribution of fire sizes and speed in our dataset.}}
\end{figure}
\begin{figure}[htbp]
    \includegraphics[width=\textwidth]{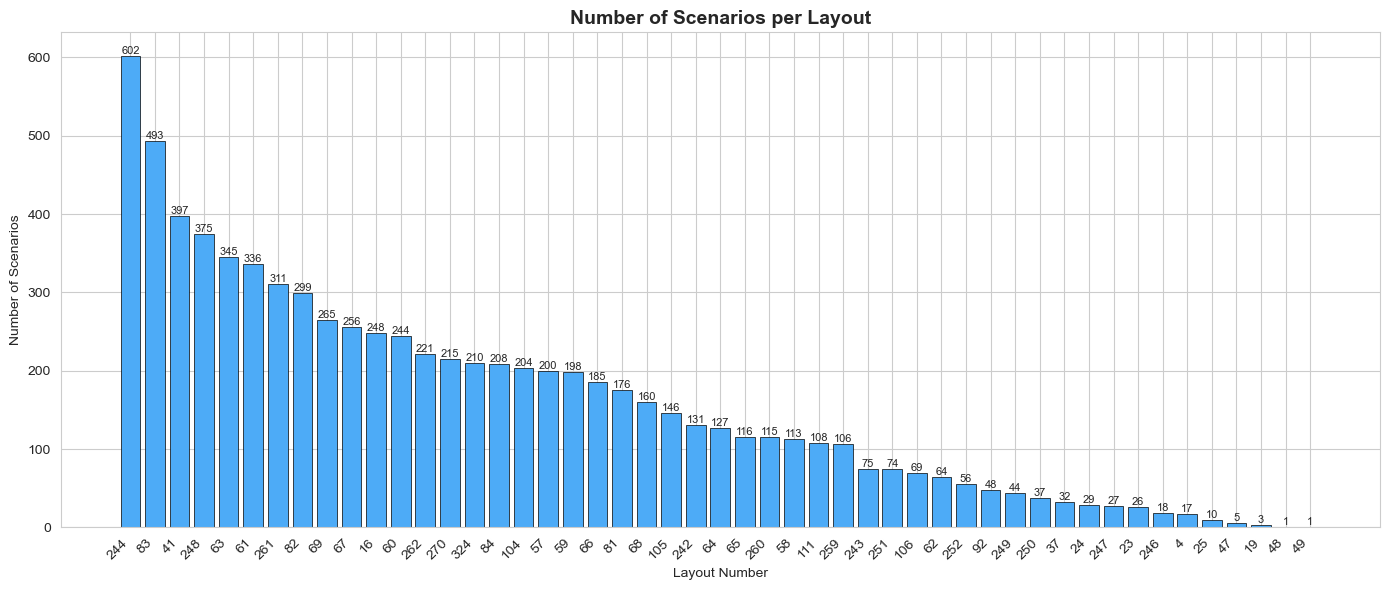}\label{fig:scenarios_per_layout}
    
    \caption{\rev{Number of wildfire scenarios per layout in our dataset.}}
\end{figure}

\section{Additional Visualizations}\label{appendix:additional visualizations}
\begin{figure}[h!]
    \centering
    
    \begin{subfigure}[b]{0.45\textwidth}
        \centering
        \includegraphics[width=\textwidth]{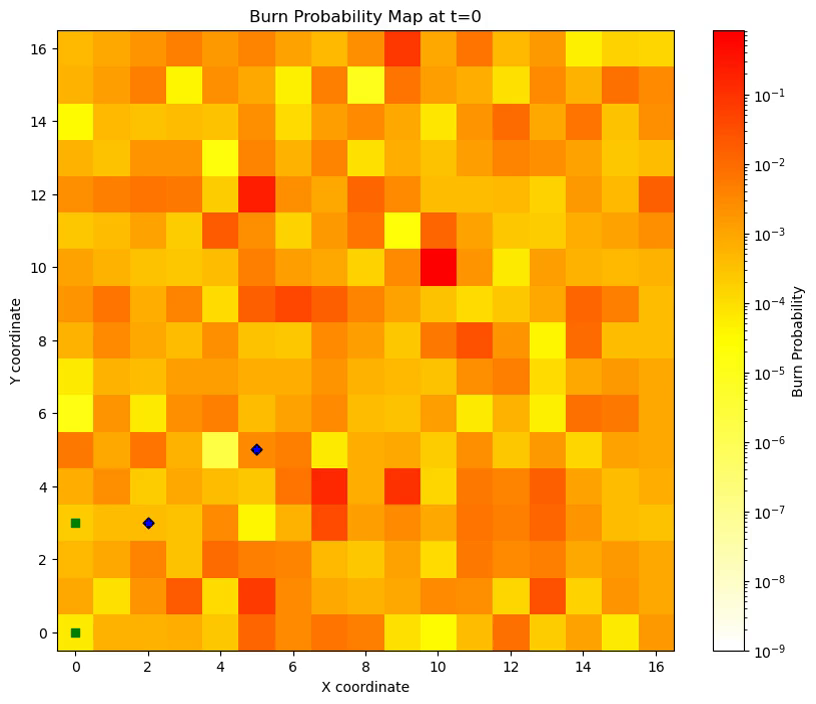}
        \label{fig:sub1}
    \end{subfigure}
    \hfill
    \begin{subfigure}[b]{0.45\textwidth}
        \centering
        \includegraphics[width=\textwidth]{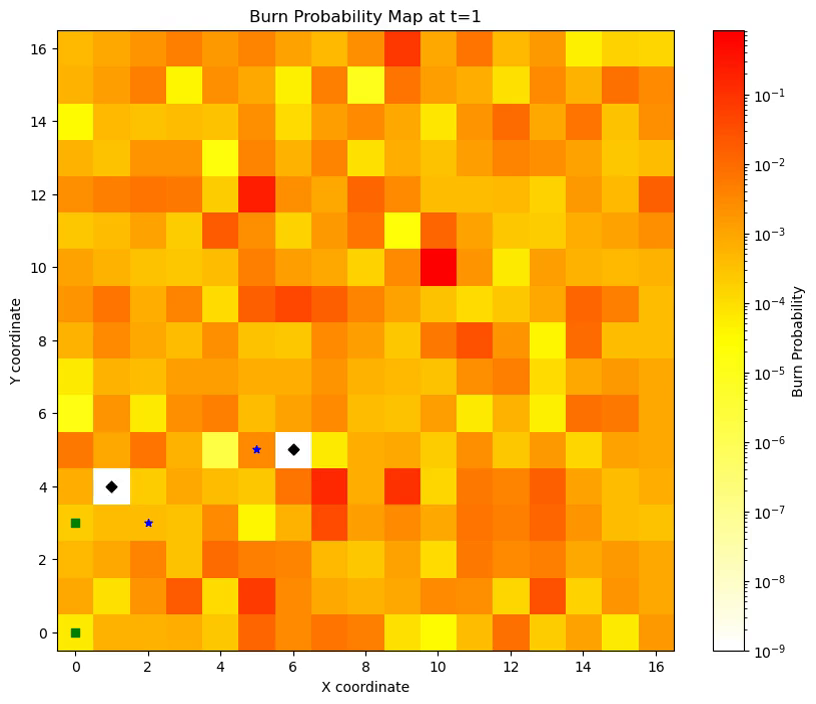}
        \label{fig:sub2}
    \end{subfigure}
    
    \vspace{0.5em}
    
    \begin{subfigure}[b]{0.45\textwidth}
        \centering
        \includegraphics[width=\textwidth]{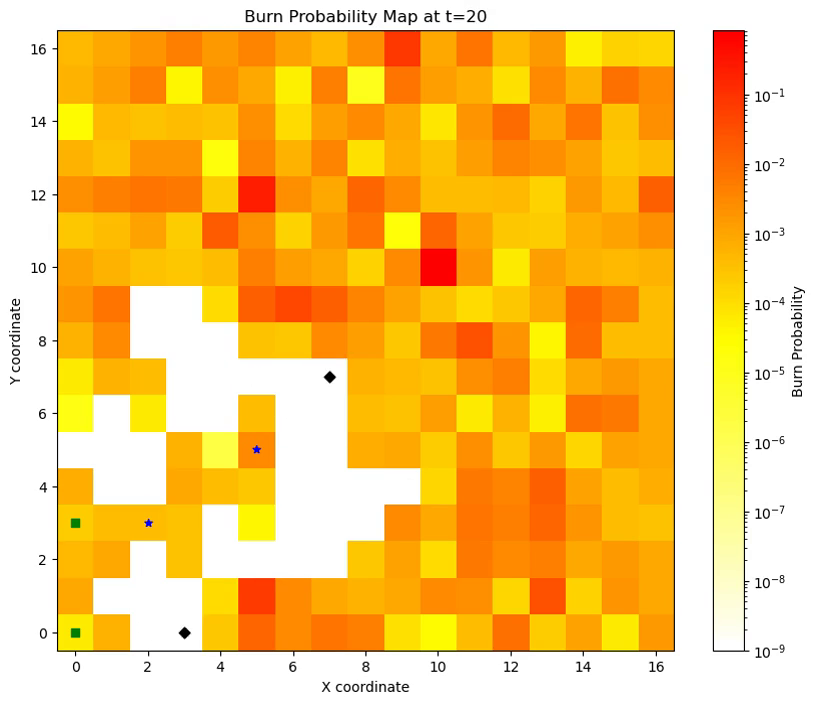}
        \label{fig:sub3}
    \end{subfigure}
    \hfill
    \begin{subfigure}[b]{0.45\textwidth}
        \centering
        \includegraphics[width=\textwidth]{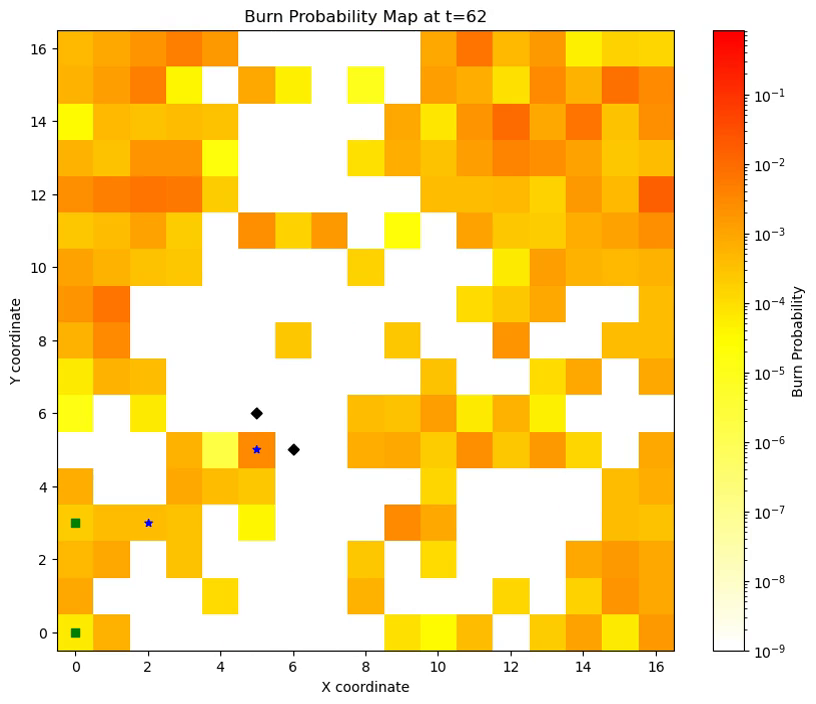}
        \label{fig:sub4}
    \end{subfigure}
    
    \caption{Visualization of ground sensors (green squares), charging stations (blue stars), and drone (black diamonds) trajectories overlaid on a small risk map. In this example, each of the two drones starts at a different charging station, but they both finish at the same one. The passage of the drones temporarily sets the visited cells' burn risk to 0. These images are frames extracted from a video generated using a function of our library.}
    \label{fig:four_images}
\end{figure}


\end{document}